\documentclass[12pt]{amsart}

\usepackage{amssymb, amscd, amsthm, verbatim, amsmath, xcolor, graphicx, turnstile, fancyhdr}
\usepackage[utf8]{inputenc} 
\usepackage[letterpaper, margin=2.5cm]{geometry}
\usepackage{lipsum} 
\title[Gevrey WKB method for P.D.O.s]{Partial Differential Operators and Fourier Integral Operators in the Gevrey setting and applications}
\author{Benjamin Colmey\textsuperscript{1} \and Richard Lascar\textsuperscript{2}}

\date{\today}

\begin{document}
\maketitle

\begin{center}
{\small
\textsuperscript{1}Department of Materials Science and Metallurgy, University of Cambridge, Cambridge, England, \textit{bc626@cam.ac.uk}\\

\medskip

\textsuperscript{2}LJAD, Université de Nice, Nice, France, \textit{richard.lascar@univ-cotedazur.fr}
}
\end{center}

\begin{abstract}
We state and prove here semiclassical results about the construction of asymptotic solutions by the WKB method for pseudo-differential equations of real principal type. It is a Gevrey version; the smooth $C^\infty$ and the analytic ones may be found in Hörmander \cite{Hormander1985} and Sjöstrand \cite{Sjoberg1982}. We present here three versions depending on the degree of entrance in the complex domain.
\end{abstract}

\renewcommand{\abstractname}{Résumé} 
\begin{abstract}
Nous présentons et prouvons ici des résultats semi-classiques concernant la construction de solutions asymptotiques par la méthode WKB pour les équations pseudo-différentielles de type principal réel. Il s'agit d'une version Gevrey ; les versions régulières $C^\infty$ et analytiques peuvent être trouvées dans Hörmander \cite{Hormander1985} et Sjöstrand \cite{Sjoberg1982}. Nous présentons ici trois versions en fonction du degré d'entrée dans le domaine complexe.
\end{abstract}

\medskip
\begin{center}
\textbf{Keywords:} Gevrey WKB method, P.D.O., real principal type\\
\textbf{Mots-clés:} Méthode WKB Gevrey, O.P.D., type principal réel
\end{center}
\medskip

\begin{center}
{\emph{Dedicated to Prof. L. Boutet de Monvel in memory.}}
\end{center}

\tableofcontents 

\section{Introduction}

In what follows, we introduce definitions of Gevrey symbols and of the related P.D.O.s. Then we state our main results.

We work in the semiclassical frame involving a small parameter $h \in (0, 1]$. We note for $m, k \in \mathbb{R},  s > 1$ by $\Sigma^{m,k}_s (\mathbb{R}^n \times \mathbb{R}^n)$ the set of Gevrey symbols, i.e., we write $a \in \Sigma^{m,k}_s$ iff for some $C > 0$ we have
\begin{equation}
|\partial_x^\alpha \partial_\xi^\beta a(x, \xi; h)| \leq C^{1 + (|\alpha| + |\beta|)} \alpha!^s \beta!^s h^{-m} \langle \xi \rangle^{k - |\beta|}
\label{eq:gevrey_symbol}
\tag{1.1}
\end{equation}
for all $\alpha, \beta \in \mathbb{N}^n$, $(x, \xi) \in \mathbb{R}^n \times \mathbb{R}^n$, $h \in (0, 1]$.

We note $\tilde{S}^m_s$ for the class $\Sigma^{m,0}_{s}$ and by $S^m_s = h^{-m} \mathcal{G}^s_b (\mathbb{R}^{2n})$

As usual, semiclassical P.D.O.s are defined by:
\begin{equation}
a(x, hD, h)u = \frac{1}{(2\pi h)^n} \int \int e^{i(x-y)\cdot \xi/h} a(x, \xi; h) u(y) \, d\xi \, dy
\tag{1.2}
\end{equation}

for a given $a \in \Sigma^{m,k}_s$. We write: $a(x, hD; h) = \mathrm{Op}_h a \quad \text{and have} \quad \mathrm{Op}_h a = \mathrm{Op}_{1}a_h, \text{ with } a_h(x, \xi) = a(x, h\xi)$.

\subsection*{Hypothesis and main results}

Let $P = P(x, hD, h)$ be a semiclassical P.D.O. of order 0 in $\mathbb{R}^n$, $n \geq 2$, i.e., with $\sigma_p(x, \xi) \in \mathcal{G}^s_b (\mathbb{R}^{2n})$. Let $p(x, \xi)$ be the principal symbol of $P$, and $H_p$ be the related Hamiltonian field. We introduce the following hypothesis near $(x_0, \xi_0) \in T^*\mathbb{R}^n$:
\begin{enumerate}
    \item[(H1)] $p$ is real,
    \item[(H2)] $dp (x_0, \xi_0) \neq 0$,
    \item[(H3)] $p(x_0, \xi_0) = 0$.
\end{enumerate}

It is customary to say that $P$ is of real principal type at $(x_0, \xi_0)$. As a heuristic, let us say that assuming more stronger assumptions than in the smooth case, we expect to find asymptotics modulo smaller remainders, i.e.,
$ \mathcal{O}(1)\exp\left(-\frac{1}{C} h^{-1/s}\right) \text{instead of } \mathcal{O}(h^\infty).$ We do obtain it, but we pay by reformulating parts of the theory.
\subsection*{Main results}

Here we shall prove 3 main theorems:

\textbf{Theorem 1.} Let $n \geq 2$, $P = p(x, hD; h)$ be a semiclassical $\mathcal{G}^s$ P.D.O. of real principal type at $(x_0, \xi_0) \in T^*\mathbb{R}^n$. Let $S \subset \mathbb{R}^n$ be a real Gevrey hypersurface $\ni x_0$ transversal to $H_p$, and let $\varphi \in G^s(\mathbb{R}^n)$ be real-valued, satisfying:
\begin{equation}
p(x, \varphi_x') = 0 \quad \text{near } x_0, \, \xi_0 = \varphi_x'(x_0).
\tag{1.4}
\end{equation}

Under these hypotheses, one may solve the WKB problem near $x_0$, i.e., if $a_0, b$ are given symbols in $S^0_s(\mathbb{R}^n)$, one may find some $a \in S^0_s(\mathbb{R}^n)$ such that:
\begin{equation}
\begin{aligned}
    &\frac{1}{h} e^{-i \varphi / h} P(a e^{i \varphi / h}) - b \in \mathcal{O}_s(h^\infty) \quad \text{close to } x_0, \\
    &a|_S = a_0.
\end{aligned}
\tag{1.5}
\end{equation}
The smallness condition in Theorem 1 is:

\paragraph{Notation 1.6.} $g(\bullet, h) \in \mathcal{G}^s_b(\mathbb{R}^n)$ is an $\mathcal{O}_s(h^\infty)$ if and only if, for some $C > 0$,
\begin{equation}
|\partial^\beta g(\bullet, h)| \leq C^{1 + |\beta|} \exp\left(-\frac{1}{C} h^{-1/s}\right)\beta!^s, \quad \forall \, \beta \in \mathbb{N}^n.
\tag{1.6}
\end{equation}

\paragraph{A convention more.} Let $A, B$ be two operators $C_0^\infty(\mathbb{R}^m)$ to $\mathcal{D'}(\mathbb{R}^n)$ defined by distribution kernels. They are said to be equivalent close to a point $(x_0, \xi_0, y_0, \eta_0) \in T^*\mathbb{R}^n \times T^*\mathbb{R}^m$ iff:
\begin{equation}
(x_0, \xi_0; y_0, \eta_0) \notin WF'_{s, h}(A - B).
\tag{1.7}
\end{equation}

In (1.7), we have used the notation $(x, \xi; y, \eta) \in WF_{s, h}(A)$ iff $(x, \xi; y, -\eta) \in WF_{s, h}(A)$, with $WF_{s, h}(u)$ being the set defined by its complement, i.e., $(x_0, \xi_0) \notin WF_{s, h}(u) $ iff
for some $\varphi \in \mathcal{G}^s_0(\mathbb{R}^n)$, $\varphi(x_0) \neq 0$, and for some neighbourhood $\mathcal{V}_{\xi_0}$ of $\xi_0$, one has:
\begin{equation}
\widehat{\varphi u}(\xi/h) = \mathcal{O}(1)\exp\left(-\frac{1}{C} h^{-1/s}\right) \text{ in } \mathcal{V}_{\xi_0}, \, \text{for } h \in (0, 1].
\tag{1.8}
\end{equation}

As for the standard WF set and the semiclassical WF set, one may use the P.D.O. calculus or FBI transforms for giving statements equivalent to (1.8).
The following result may be called the Gevrey Egorov Theorem. The statement is:

\textbf{Theorem 2.} Let $P = P(x, hD; h)$ be a $\mathcal{G}^s$ P.D.O. of order 0 on $\mathbb{R}^n$, $n \geq 2$. Assume $P$ has a real principal symbol $p(x, \xi)$ satisfying at $(x_0, \xi_0) \in T^*\mathbb{R}^n$, $p(x_0, \xi_0) = 0$, $H_p(x_0, \xi_0) \neq 0$. Then, there is a $G^s$ canonical map $\chi$ from a neighbourhood of $(0,0)$ to a neighbourhood of $(x_0, \xi_0)$ such that:
\begin{equation} \tag{1.9}
\begin{aligned}
&\text{(i)} \quad \chi^* p = \xi_1, \\
&\text{(ii)} \quad \text{There are } \mathcal{G}^s \text{ F.I.O. } A \text{ and } B \text{ related to } \chi \text{ and } \chi^{-1} \text{ such that:} \\
&\quad \text{(a)} \quad WF'_{s, h}(A) (\text{ resp } WF'_{s, h}(B)) \text{ is close to } (x_0, \xi_0; 0, 0) \text{ (resp.} (0, 0, x_0, \xi_0)\text{)}, \\
&\quad \text{(b)} \quad (x_0, \xi_0) \notin WF'_{s, h}(AB - I) \text{ and } (0, 0) \notin WF'_{s, h}(BA - I), \\
&\quad \text{(c)} \quad (x_0, \xi_0; x_0, \xi_0) \notin WF'_{s, h}(h A D_{1} B - P) \text{ and } \\
&\quad \qquad (0, 0; 0, 0) \notin WF'_{s, h}(B P A - h D_{1}).
\end{aligned}
\end{equation}

In other words, $P$ and $hD_{1}$ are $\mathcal{G}^s$ microlocally conjugate near $(x_0, \xi_0)$.

Before stating Theorem 3, we outline FBI transforms in the spirit of Sjöstrand \cite{Sjoberg1982}.

Let $x_0 \in \mathbb{C}^n$, $y_0 \in \mathbb{R}^n$, $\eta_0 \in \mathbb{R}^n \setminus 0 $, be given, $\varphi(x, y)$ an holomorphic map close to $(x_0, y_0)$ with:
\begin{equation}
\varphi'_y(x_0, y_0) = -\eta_0 \in \mathbb{R}^n, \quad \mathrm{Im} \, \varphi''_{y,y}(x_0, y_0) > 0, \quad \det \varphi''_{xy}(x_0, y_0) \neq 0.
\tag{1.10}
\end{equation}

FBI transforms are integral transforms mapping distributions close to $y_0$ to holomorphic functions close to $x_0$ of the type:
\begin{equation}
T u(x, h) = \int e^{\frac{i}{h} \varphi(x, y)} a(x, y; h) u(y) \chi(y)\, dy,
\tag{1.11}
\end{equation}
where in (1.11), $a(x, y; h)$ is an analytic symbol elliptic at $(x_0, y_0)$, and $\chi \in C_0^{\infty}(\mathbb{R}^n)$ is 1 close to $y_0$.
Following \cite{Sjoberg1982}, we set:
\begin{equation}
\phi(x) = sup \underset{\substack{y \text{ neighbourhood of } \\ (y_0, \mathbb{R}^n)}} { } \mathrm{-Im} \, \varphi(x, y).
\tag{1.12}
\end{equation}

$\phi$ is smooth, plurisubharmonic, and $\Lambda_\phi = \{ (x,\frac{2}{i} \partial_x \phi(x); x \in \text{ 
 neighboorhood } (x_0, \mathbb{C})\} $ is I Lagrangian and $\mathbb{R}$-symplectic. The complex canonical transform related to $T$ is:
\begin{equation}
K_T : (y, -\varphi'_y (x, y)) \mapsto (x, \varphi_x'(x, y)).
\tag{1.13}
\end{equation}

$K_T|_{\mathbb{R}^{2n}}$ maps $\mathbb{R}^{2n}$ to $\Lambda_\phi$.

The basic example is the Bargmann transform $T = T_0$, $\varphi = \varphi_0 = \frac{i}{2}(x - y)^2$, $\phi_0 = \frac{1}{2}(Im (x))^2$, a $\equiv$ 1, and $K_{T_0}|_{\mathbb{R}^{2n}} : (y, \eta) \to (y - i\eta, \eta)$.
\begin{equation}
\tag{1.14}
\end{equation}

\paragraph{A definition more.} Let $x_0 \in \mathbb{C}^n$, $F(x, h)$ a smooth map close to $x_0$. We note $F \sim 0$ in $\mathcal{G}^{s'}$ for the weight $\phi$ if, in a neighbourhood $\mathcal{V}$ of $x_0$ $\in \mathbb{C}^n$, one has for some $C > 0$:
\begin{equation}
|\partial^\alpha F(x, h)| \leq C_\alpha e^{\frac{1}{h}{\phi}(x)} \exp\left(-\frac{1}{C} h^{-1/s'}\right),  \quad h \in (0, 1].
\tag{1.15}
\end{equation}

\textbf{Theorem 3.} Let $P(y, D_y)$ be a $\mathcal{G}^s$ differential operator of degree $m$ such that near $(y_0, \eta_0) \in T^*(\mathbb{R}^n \setminus 0)$, $p(y, \eta)$ the principal symbol of $P$, is of real principal type. We have:

(i) There is a Gevrey $s$ FBI transform defined near $y_0 \in \mathbb{R}^n$, i.e., an integral transform of type (1.11) with some phase function $\varphi(x, y)$ of $G^s$ class close to $(x_0, y_0)$, with $x_0 = (y_0 - i(0, \eta'_0) \in \mathbb{C}^n$, and where:
\begin{equation}
\varphi'_y(x_0, y_0) = -\eta_0, \quad \mathrm{Im} \, \varphi''_{y,y}(x_0, y_0) > 0, \quad \det \varphi''_{xy}(x_0, y_0) \neq 0,
\tag{1.16}
\end{equation}

a symbol $a(x, y; h)$ defined in $\mathbb{C}^n \times \mathbb{R}^n$ close to $(x_0, y_0)$ of $G^s$ class elliptic at $(x_0, y_0)$ such that, if $u \in \mathcal{D'}(\mathbb{R}^n)$ is independent of h $>$ 0, one has: 
\begin{equation}
    (h D_{Re x_1} T_{u} - h^m T P(y, D_y)u) (x, h) \sim 0 \tag{1.17}
\end{equation}
for the weight $\phi_0 = \frac{1}{2}(\mathrm{Im} \, x)^2$ in $\mathcal{G}^{2s-1}$ close to $x_0 \in \mathbb{C}^n$.

(ii) The map $K_T : (y, -\varphi'_y) \mapsto (x, \varphi'_x)$, restricted to $T^*\mathbb{R}^n \setminus 0$, is a canonical transform $\chi$ from a neighbourhood of $(y_0, \eta_0)$ to a neighbourhood of $(x_0, \xi_0)$, $\xi_0 = \frac{2}{i} \partial_x \phi_0(x_0)$ in $\Lambda_\phi$, the I-Lagrangian space related to the Bargmann transform, such that:
\begin{equation}
\xi_1|_{\Lambda_\phi} = p_0 \chi^{-1}
\tag{1.18}
\end{equation}

In few words, Theorem 3 above is an attempt of intertwining $h^mP(y,D_y)$ to $h D_{Re(x_1)}$ by an integral transform of FBI type. In the travel into the complex domain, we have had to take almost holomorphic extensions with respect to:
\begin{equation}
\Gamma = \{(x, y(x)) \in \mathbb{C}^n \times \mathbb{R}^n; x \in \text{neigh}(x_0, \mathbb{C}^n)\} \text{ to } \mathbb{C}^n \times \mathbb{C}^n.
\tag{1.19}
\end{equation}

Having defined $y(x)$ by $y(x) = \pi_x (\chi^{-1} (x,\xi(x))$, with $\pi$ the canonical projection $T^*(\mathbb{R}^n) \to \mathbb{R}^n$, and $\xi(x) = \frac{2}{i} \partial_x \phi_0(x)$.

For us, this explains the loss $2s-1$ in the Gevrey setting of (1.17).
\subsection*{Acknowledgments}
The second author is very grateful to Johannes Sjöstrand for helpful discussions.
\section{Stationary phase and symbol calculus in the Gevrey setting}

We look at the smallness of phase integrals of type:
\begin{equation}
J_f(y, h) = \int_{\mathbb{R}^n} e^{\frac{i}{h} f(x, y)} a(x, y; h) \, dx.
\tag{2.1}
\end{equation}

In (2.1), it is assumed $f$ is complex-valued of $G^s$ class on $\mathbb{R}^n \times \mathbb{R}^m$, i.e., $f \in S^0_s(\mathbb{R}^n \times \mathbb{R}^m)$, with $\mathrm{Im} \, f(x, y) \geq 0$, and that $a \in S^{m_0}_s(\mathbb{R}^n \times \mathbb{R}^m)$ has $x$-compact support.

One has the lemma: 
\subsection*{{Lemma 2.1}} (Non-stationary phase lemma)
Assume $\mathrm{Im} \, f \geq 0$ and 
\begin{equation}
\mathrm{Re} \, f'_x(x, y) \neq 0 
\quad \text{for } (x, y) \in \mathrm{supp} \, a 
\quad \text{when } \mathrm{Im} \, f = 0.
\tag{2.2}
\end{equation}

Then, $J_f(y, h)$ is a small $\mathcal{G}^s$ remainder, i.e., one has locally uniformly in $y$ estimates, for some $C > 0$, all $\beta \in \mathbb{N}^m$, and $a \in (0, 1]$:

\begin{equation}
|\partial^\beta _y J_f(y, h)| \leq C^{1+ |\beta|} \beta!^s  \exp\left(-\frac{1}{C} h^{-1/s}\right).
\tag{2.3}
\end{equation}

\paragraph{Proof.} We perform an IPP in (2.1), referring to \cite{LascarLascar1997} for a justification in the Gevrey case. Let us sketch briefly the argument. Using hypothesis (2.2), we introduce the IPP operator $L$:
\begin{equation}
L(a) = \sum_{j=1}^n \frac{1}{i} h \partial_{x_j} \left(\frac{\bar{f'_{x_j}}}{| f_x'|^2} a\right)
\tag{2.4}
\end{equation}

Using the Gevrey quasinorms of Lemma 1 in \cite{LascarLascar1997} (with the standard metric $|t|_{x} = |t| \, \forall \, x$), we write for  $N > 0$ and for $\epsilon > 0$ small

\[
N(L^N(a), T)(x, y) \leq 
\left[\sum_{j=1}^n N(f_{j}, T\left(1 + \frac{\epsilon}{N}\right)^{sN})(x, y)\right]^{N}
N(a, T(1 + \frac{\epsilon}{N})^{sN}(x, y) (T^{-N}\epsilon^{-sN}N^{sN})
\tag{2.5}
\]
with $f_j = h \bar{f}_{x_j}' / |f'_x|^2$. Since one has:
\[
N(f_j, T)(x, y) \leq h |f'_x|^{-2} N(\bar{f}_{x_j}', T) (1 - \dot{N}(|f'_x|^2,T)(|f'_x|^{-2})^{-1} \quad \tag{2.6}
\]

We use assumption (2.2) and in (2.6), $T \leq C |f_{x}'|$, and we bound:$N(|\bar{f}_{x_j}'|, T) \leq C |f_{x}'|$,
and we note that the third term in (2.6) is bounded so that we have a factor: 
\begin{equation}
    N^{sN} h^N |f_{x}'|^{-N} |f_{x}'|^{-N} T_0^{-N} \text{   For some small   } T_0 > 0
    \tag{2.7}
\end{equation}

So performing $N$ integration by parts, we have a bound:$N^{sN} h^N (T_0 / C^2)^{-N}$ and taking $N = N_\text{min}$, we are done with (2.3). We also have the lemma:

\subsection*{Lemma 2.2} (Gevrey stationary phase lemma in the real case)
We assume at $(x_0, y_0) \in \mathbb{R}^n \times \mathbb{R}^m$ the following hypotheses:
\begin{equation}
\begin{aligned}
&\text{(a)} \quad f \text{ is real-valued close to } (x_0, y_0), \\
&\text{(b)} \quad f_x'(x_0, y_0) = 0, \\
&\text{(c)} \quad \det f_{x,x}''(x_0, y_0) \neq 0.
\end{aligned}
\tag{2.8}
\end{equation}

(a), (b), and (c) imply the existence of a $G^s$ map $y \mapsto \chi(y)$ close $y_0$ such that:
\begin{equation}
f_x'(\chi(y), y) = 0, \quad \chi(y_0) = x_0. \tag{2.9}
\end{equation}
And one has:
\begin{equation}
b(y, h) = e^{\frac{-i}{h} f(\chi(y), y)} \int a(x, y; h) e^{\frac{i}{h} f(x, y)} \, dx, \tag{2.10}
\end{equation}
is a $G^s$ symbol in $\mathbb{R}^m$ close to $y_0$, of order $m_0 - n/2$, enjoying the asymptotic expansion.

\begin{equation}
b(y, h) \sim \left|\det \frac{f_{x,x}''(\chi(y), y)}{2\pi h}\right|^{-\frac{1}{2}} 
e^{i\frac{\pi}{4} \sigma} \sum_{j=0}^\infty h^j L_{f,j,y} a(\chi(y), y), \tag{2.11}
\end{equation}
where $L_{f,j,y}$ is a differential operator of degree $\leq 2j$ in $x$, and $\sigma = sgn(f''_{x,x} (x_0, y_0)$.

Moreover, the asymptotic expansion in (2.11) above in the RHS is only a formal symbol of Gevrey type $2s-1$, as pointed out in earlier works.

\paragraph{Proof.} Since $f$ is real-valued, we may use the Morse Lemma; By Taylor expansion at order 2 close $(x_0, y_0)$, we write:
\[
f(x, y) = f(\chi(y), y) + \frac{1}{2} < Q(x, y)(x - \chi(y)), x - \chi(y) >
\]
where:
\[
Q(x, y) = 2 \int_0^1 (1 - t) \, f_{x,x}''(\chi(y) + t(x - \chi(y)) y) \, dt.
\]

Setting $Q_0 = f''_{x,x}(x_0, y_0)$, one has by the Morse lemma a map $(x, y) \longrightarrow R(x, y)$ such that:
\[
Q_0 = \quad ^{t}R(x, y) Q(x, y) {R}(x, y), \quad {R}(x_0, y_0) = \mathrm{Id}.
\]

Performing the change of variables:
\[
\tilde{x}(x, y) = R^{-1}(x, y)(x - \chi(y)), \quad \tilde{y} = y,
\]
near $(x_0, y_0)$, we write:
\[
b(y, h) = \int e^{\frac{i}{h} (Q_0 \tilde{x},\tilde{x})}\tilde{a}(\tilde{x}, y, h) \, d\tilde{x},
\]
for a symbol $\tilde{a}$ in $G^s$ with an $\tilde{x}$ support close to 0. Next, using the Plancherel's formula, one has:
\[
b(y, h) = (\det(Q_0/2i\pi h))^{-1/2} e^{\frac{-i h}{2} (Q^{-1}_0D,D)} \tilde{a}(0, y).
\]

Now the asymptotic expansion (2.7) follows from Hörmander \cite{Hormander1985}, Theorem 7.6.9. When $f$ is complex-valued, the situation is more intricate. We refer to \cite{Hormander1985} and \cite{HelmSjoberg1975} for the smooth case, and to \cite{LascarLascar1997} for the Gevrey one. Here, in the proof of Theorem 3, we shall adapt the method of \cite{HelmSjoberg1975}. Now we can go to the PDO symbolic calculus.

{Symbolic calculus in the Gevrey setting}

We shall treat here the Gevrey classes $\Sigma^{m, k}_s(\mathbb{R}^n \times \mathbb{R}^n)$ as it is known for the classes $S^m_s(\mathbb{R}^n \times \mathbb{R}^n)$ see \cite{Lascar1988}, \cite{Hitrik2023complex}.  

One has to check first that if $q \in \Sigma^{m_0, k_0}_s$, $a \in \Sigma^{m', k'}_s$, the composition:
\[
(q \circ a)(x, \xi) = \frac{1}{(2\pi h)^n} \int\int e^{-i (x - y) (\xi - \eta) / h} q(x, \eta) a(y, \xi) \, dy \, d\eta \tag{2.12}
\]
is in the class:\[ \Sigma^{m_0 + m', k_0 + k'}_s(\mathbb{R}^n \times \mathbb{R}^n). \]

\subsection*{{Proposition 2.3}}
Let $Q = \mathrm{op}_hq$, $A = \mathrm{op}_ha$, where $q \in \Sigma^{m_0, k_0}_s$, $a \in \Sigma^{m', k'}_s$, given then $Q \circ A = \mathrm{op}_h(q \circ a)$ with $(q \circ a)$ given by (2.12). One has:$(q \circ a) \in \Sigma^{m_0 + m', k_0 + k'}_s$
and for every $N > 0$:
\[
(q \circ a)(x, \xi) = \sum_{|\alpha| < N} \frac{h^{|\alpha|}}{\alpha!} D_\xi^\alpha q(x, \xi) \partial_{x}^\alpha a(x, \xi) + r_N(q, a)(x, \xi), \tag{2.13}
\]
with $r_N(q, a) \in \Sigma^{m_0 + m' - N, k_0 + k' - N}_s$, given by:
\[
r_N(q, a)(x, \xi) = h^N \sum_{|\alpha| = N} \frac{N!}{\alpha!} \int\int_{[0, 1]\times{R^{2n}}} \frac{(1 - s)^{N-1}}{(2\pi h)^n(N-1)!} 
e^{-i \frac{y\eta}{h}} D_\xi^\alpha q(x, \xi + \eta) \partial_x^\alpha a(x + sy, \xi) \, ds \, dy \, d\eta \tag{2.14}
\]
\paragraph{Proof.}
The proof relies on the following definition and follows the lines of \cite{Lerner2011} Lemma 4.12 and 4.15 and from a partition of unity in $\xi$ variables in the spirit of Hörmander \cite{Hormander1985}.

\subsection*{{Lemma 2.4}}
For $t \in \mathbb{R} \setminus \{0\}$, we define the map $J^t$ by
\[
J^t b = \frac{1}{(2\pi |t|)^n} \int\int_{\mathbb{R}^{2m}} e^{-i z \zeta / t} b(x + z, \xi + \zeta) \, dz \, d\zeta, \quad b \in \Sigma^{m', k'}_s(\mathbb{R}^n \text{x} \mathbb{R}^n). \tag{2.15}
\]
Then $J^t$ maps $\Sigma^{m', k'}_s$ into itself, and for any $N > 0$, one has:

\[
J^t b(x, \xi) = \sum_{|\alpha| < N} \frac{t^{|\alpha|}}{\alpha!} D^\alpha_\xi D^\alpha_x b(x, \xi) + r_N(t)(x, \xi), \tag{2.16}
\]
where \( r_N(t)(x, \xi) \in \Sigma^{m'-N, k'-N}_s \), and 
\[
r_N(t) = t^N \int_0^1 \frac{(1 - s)^{N-1}}{(N-1)!} J^{st} \left((D_\xi \partial_x)^N b \right)(x, \xi) ds. \tag{2.17}
\]

\paragraph{Proof.} The proof uses a dyadic partition of unity and other results of Hörmander \cite{Hormander1985} section 7.7. 

First, observe that \( J^h \) is a Fourier multiplier. 
\[
J^h b(x, \xi) = (2\pi h)^{-n} \iint_{\mathbb{R}^{2n}} e^{-i y \cdot \eta / h} b(x + y, \xi + \eta) \, dy \, d\eta,
\]
and hence 
\[
J^h b = (2\pi h)^{-n} \left[ \left(e^{\frac{-i}{2h} (B\cdot,\cdot)}\right) * b \right](x, \xi). \tag{2.18}
\]

\[
B = 
\begin{pmatrix}
0 & I_n \\
I_n & 0
\end{pmatrix}, \quad (B X, X) = 2 x \cdot \xi \quad \text{if } X = (x, \xi).
\]
One has \( B^{-1} = B \), \(\text{sgn}(B) = 0\), and \(\det(B) = 1\). Hence,
\[
J^h b = e^{\frac{ih}{2} (B D, D)} b \quad \text{since }  \mathcal{F}^{-1}_{ \Xi  \rightarrow X} e^{(\frac{ih}{2}) (B\Xi,\Xi  )} = (2 \pi h)^{-n} e^{-(i/2h ) (B X, X)}.
\] which is a standard fact.

From (2.18), it can be seen that \( J^h \) commutes with derivations and integration by parts up to the dimension gives
\[
\left| \partial^\gamma J^h b(x, \xi) \right| \leq C_n (1 + h)^{n+2} \sup_{|\alpha|, |\beta| \leq n+2} | \partial^\alpha_x \partial^\beta_\xi \partial^\gamma b |_{L^\infty}, \tag{2.19}
\]
for some constant \( C_n \) depending only on $n \geq 1$

A consequence is that \( J^h \) maps continuously \( C^\infty_b \) into itself and \( \mathcal{G}^s_b \) unto itself (this last point uses Gevrey quasi-norms).

Let us consider now a symbol \( b \in \Sigma_s^{0, k}, \, k \in \mathbb{R} \), i.e., satisfying
\begin{equation}
|\partial_x^\alpha \partial_\xi^\beta b(x, \xi)| \leq C^{ (1 + |\alpha| + |\beta|)}  \alpha^s! \beta^s!\langle \xi \rangle^{k - |\beta|} \tag{2.20}
\end{equation}

We want to check that \( J^h b \) is in the same class. First, we cut the integral \( J^t \) close to the critical point \( y_c = \eta_c = 0 \) and write \[J^h b = I_\chi b + r_\chi b\] having split \( 1 = \chi_r(\Upsilon) + (1 - \chi_r )(\Upsilon) \). For some \(\chi_r\) with support in \( |\Upsilon| \leq r \) in \( |\Upsilon|\leq r/2 \), we evaluate each term by the stationary and non-stationary phase lemma. Since the parameter ($\xi$) is in \(\mathbb{R}^n\), we introduce the dyadic partition of unity of \cite{Hormander1985}. There are \( \mathcal{G}^s \) functions with compact supports \(\chi_0 \in \mathcal{G}_0^s(\mathbb{R}^n)\), \(\chi \in \mathcal{G}_0^s(\mathbb{R}^n\setminus 0)\)such that
\[
1 = \chi_0 + \sum_{\nu \geq 1} \chi_\nu,\tag{2.21}
\]
\[
\chi_\nu(\cdot) = \chi(2^{-\nu}\cdot)
\]

For \( \nu = 0 \), we have
\[
\chi_0 I^h b = I_h^0 + r_\chi^0, \tag{2.22}
\]
and for \(\nu\geq 1 \), we perform a change of variables \( \tilde{\eta} = \eta/R_\nu, \, R_\nu = 2^\nu\) and set
\[
\chi_\nu J_h^\ell b = I^{\nu}_{\tilde{h}} + r_{\tilde{h}}^{\nu},  \tag{2.23}
\]
having set
\[
\tilde{h} = \frac{h}{R_{\nu}}\leq 1,
\]
and
\[
I_{\tilde{h}}^{\nu} = \chi_{\nu}(\xi) \int \int e^{-i y\tilde{\eta}/ \tilde{h}} b(x + y, \xi + R_\nu \tilde{\eta}) \chi(y, \tilde{\eta}) \, \tilde{\eta} \frac{dy d\tilde{\eta}}{(2\pi \tilde{h})^n}
\]
\[
r_{\tilde{h}}^{\nu} = \chi_{\nu}(\xi) \int \int e^{-i y\tilde{\eta}/ \tilde{h}} b(x + y, \xi + R_\nu \tilde{\eta}) (1-\chi)(y, \tilde{\eta}) \, \tilde{\eta} \frac{dy d\tilde{\eta}}{(2\pi \tilde{h})^n}, \tag{2.24}
\]

Hence, we have
\[
J^h b = \sum_{\nu \neq 0} \chi_\nu J^h(b) + \chi_0 J^h b,
\]
which is a sum of terms.

The symbol \( I_h^o(x, \xi) \in S_s^{o}(\mathbb{R}^n \times \mathbb{R}^n) \) but has compact support in \(\xi\). Hence,
\[
I_h^{o} (x, \xi)\in \Sigma_s^{o, k}(\mathbb{R}^n \times \mathbb{R}^n), \quad r_\chi^{o} \text{ is a small } \mathcal{G}^s \text{ remainder by the non-stationary phase lemma.}
\]

We address the reader to \cite{Hitrik2023complex} for the estimate \((3.149)\), which gives
\[| \partial^\alpha {r^o_{\chi}}| \leq C^{ (1 + |\alpha|)} \alpha!^s e^{-\frac{1}{C} h^{-1/s}} \langle \xi \rangle^{-|\alpha_\xi| +k}. \tag{2.25}\]
since $<\xi> \sim \text{on supp} \chi_0$.

In relation \((2.24)\), we may apply the stationary phase lemma with the quadratic phase \(y \cdot \tilde{\eta}\) and small parameter $\tilde{h} = h / R_{\nu}$, and we deduce from \cite{Hitrik2023complex}, relations \((3.150)\) and \((3.151)\), if $I_{\tilde{h}}^\nu$ is defined in \((2.24)\), if we set \( \tilde{\xi} = \xi / R_\nu \), for integers \(K > 0 \),
\[
\tilde{R}_{K, \tilde{h}}(x, \tilde{\xi}) = {I^\nu}_{\tilde{h}}(x, \tilde{\xi}) - C_B \sum_{k=0}^{K-1} \frac{\tilde{h}^k}{k! (2i)^k} ({B}^{-1} \tilde{D},\tilde{D})^k \tilde{b}_\nu| y = \tilde{\eta} = 0,
\]
with
\[
\tilde{b}_\nu = b(x + y, R_\nu \tilde{\xi} + R_\nu \tilde{\eta}) \chi(y, \tilde{\eta}),\quad C_B = e^{\frac{i \pi}{4} \operatorname{sgn}B} |\det B|^{-1/2}.\]

We check
\[
|\partial^\gamma_{x,\tilde{\xi}} \tilde{R}^\nu_{K, \tilde{h}}(x, \tilde{\xi})| \leq C^{(1 + |\alpha| + K)} K!^{2s-1}!\gamma^s!  \tilde{h}^K
\tag{2.26}
\]
Indeed, applying \cite{Hitrik2023complex} inequality (3.151), one has to check that \( R_\nu^{-k} \tilde{b}_{\nu}(x, y, \tilde{\xi},  \tilde{\eta}) \) is bounded in \( \mathcal{G}_b^s \big( \mathbb{R}_{x,\tilde{\xi}}^{2n} \times \mathbb{R}_{y,\tilde{\eta}}^{2n} \big) \),  which is a consequence of the assumption: \( b \in \tilde{S}_s^{0, k} (\mathbb{R}^{n} \times \mathbb{R}^{n}) \). Since when taking derivatives w.r.t. \( (\tilde{\xi}, \tilde{\eta}) \), we obtain a factor \( R_\nu^{|\gamma_{\tilde{\xi}}| + |\delta \tilde{\eta}|} \),  which is canceled by the factor \( \langle R_\nu \tilde{\xi} + R_\nu \tilde{\eta} \rangle^{k - |\gamma_{\tilde{\xi}}| - |\delta \tilde{\eta}| } \). Since \( |\tilde{\xi}| \sim 1 \), \( |\tilde{\eta}| \leq r \) small. s.t. \( |{\xi} + \eta| \sim |\xi|\).

We deduce for \[R^\nu_{K, \tilde{h}}(x, \xi) = \chi (\xi / R_\nu) \tilde{R^\nu}_{{K}, \tilde{h}}(x, \xi/R_\nu),
\]
estimates
\[|\partial_{x,\xi}^\gamma {R^\nu}_{K, \tilde{h}}(x, \tilde{\xi})| \leq C^{(1 + |\gamma| + K)} K^{2s-1}!\gamma^s! |\xi|^{-K -|\gamma_\xi|} h^{K}|\xi|^{k}.
\tag{2.27}\]

(use \( R_\nu \approx |{\xi}| \) on \(\mathrm{supp} \, \chi_\nu\) and Leibniz formula) uniformly in $\nu$ and $K$.

We also have, by applying \cite{Hitrik2023complex},
\[|\partial^\gamma r^\nu_{\tilde{h}}| \leq C^{1+|\gamma|} \gamma!^{s} \exp \left( -\frac{1}{C} \left( \frac{|{\xi}|}{h} \right)^{1/s} \right) \ \text{uniformly in } \nu. \tag{2.28}
\]

Observe now that from (2.27) that the sum \[R_{K, {h}} = \sum_{\nu \neq 0}^\infty {R^\nu}_{K, \tilde{h}}\]
is well-defined in \(\tilde{S}_s^{-K, -K + 1/2}\) for \( K \geq 1 \), since \( |\xi| \sim R_\nu \sim 2^\nu \) on \(\mathrm{supp} \, \chi_\nu\).

Now, using asymptotics (2.27) and \cite{Hitrik2023complex} (3.151) for \( I^\nu_{\tilde{h}} \), one has, changing \( K \) into \( K + 1 \), 
\[
|\partial^\gamma \left( J^h b - \sum_{\nu = 0}^{K-1} h^k / k!  (D_\xi \partial_x)^k b(x,\xi) \right) | 
\leq C^{1 + |\gamma| + |K|} {h^K}{K!}^{2s-1} \gamma!^s \langle \xi \rangle^{-K - |\gamma_\xi|+k} \tag{2.29},
\]
since we may write \[R^K_{h} = R^{K+1}_{h} + \frac{h^K}{K!} (D_\xi \partial_x)^K b, \]
which is the desired result. So Lemma (2.4) is proved. We can now go to the end of Theorem 2's proof.

\section{Proof of Theorem 2}

Let \( a(x, \xi', h) \in \tilde{S}_s^{m, k}(\mathbb{R}^n \times \mathbb{R}^{n-1}) \) be a PDO in \(\mathbb{R}^n\) given by \( a \), i.e.,
\[
a(x, h D_x,h) u = (2 \pi h)^{-(n-1)} \int e^{i x' \xi' / h} a(x, \xi', h) \hat{u}(x_1,\xi'/h) \, d \xi' \tag{3.1}.
\]

We introduce some classes of Gevrey FIO in setting:
\[Fu (x, h) = (2 \pi h)^{-(n-1)} \int e^{\frac{i}{h} S(x,\eta')} a(x, \eta', h) \hat{u}(x_1,\eta'/h) \, d \eta' \tag{3.2},\]
with symbols \( a \in \tilde{S}_s^{m', k'}(\mathbb{R}^n \times \mathbb{R}^{n-1}) \) and phase given by \(S(x, \eta') \in S_s^{0, 1}(\mathbb{R}^n \times \mathbb{R}^{n-1})\), with $S$ real and 
$det S''_{x' ,\eta'} \neq 0 \, \text{(see Eskin)}$ \cite{Eskin2011} for the classical homogeneous case).
Let \( P(x, h D_x, h) \) be a \(\mathcal{G}^s\) PDO on \(\mathbb{R}^n\) of bi-order \((0,1)\) with \( p(x, \xi) \) as a principal symbol of real principal type at \((x_0, \xi_0) \in T^* \mathbb{R}^n \setminus 0 \). More precisely, we assume \( H_p (x, \xi_0) \) transversal to the fiber at \((x_0, \xi_0)\). Using the implicit function theorem, we may argue:

\[
P = h D_{x_1} + Q(x, h D_{x!}, h) \tag{3.3}
\]

\( Q \) being a \(\mathcal{G}^s\) PDO of bi-order \((0,1)\) with principal symbol $-\lambda(x, \xi'), \quad \lambda(x, \xi') \in \tilde{S}_s^{0,1}$, $\text{real, so } \, p(x, \xi) \, \text{rewrites}$

\[
p(x, \xi) = \xi_1 - \lambda(x, \xi'). \tag{3.4}
\]

We solve close to \((x_0, \xi_0')\) the eikonal equation 

\[
\begin{cases}
    \frac{\partial \varphi}{\partial x_1} = \lambda(x, \varphi'_{x'}), \\
    \varphi|_{x_1 = 0} = x' \xi',
\end{cases} \tag{3.5}
\]

In view of proving the Egorov Theorem, one has to use FIO. We have written \((x_0, \xi_0) = (x_0, \lambda(x_0, \xi_0'), \xi')\), \(\varphi\) is \(\mathcal{G}^s\) and \underline{real} . Let \(F\) be the h FIO:

\[
F u(x, h) = (2 \pi h)^{-(n-1)} \int e^{i \frac{\varphi(x, \eta')}{h}} \alpha(x, \eta'; h) \hat{u}(x, \eta'/h) \, d\eta', \tag{3.6}
\]

We assume \( a(x, \eta') \in \tilde{S}_s^{m, k} \) elliptic at \((x_0, \xi_0')\) and supported in a neighbourhood of this point.

We compute now
\[
(P(x, h D_x, h) F u = h F D_{x_1} u + F_1 u \tag{3.8}
\]
where \( F_1 \) is similar to \( F \), its symbol is \( a_1(x, \eta') \), given by
\[
a_1(x, \eta') = e^{\frac{-i}{h} \varphi(x, \eta')} Q(a e^{\frac{i}{h} \varphi}) + \frac{h}{i} \partial_{x_1} a + \varphi'_{x_1} a, \tag{3.9}
\]
One has \( a_1 \in \tilde{S}_s^{m-1, k} \) since we solved in (3.5) the eikonal equation. Indeed, it has the expansion 
\[
\sum_{\alpha \in N^n} \frac{h^{|\alpha|}}{\alpha!} \sigma_p^{(\alpha)}(x, \varphi'_{x'}) D_y^\alpha \left(a e^{i\frac{\rho}{h}} \right) \big|_{y = x}, \tag{3.10}
\]
with
\[
\rho(y, x, \eta') = \varphi(y, \eta') - \varphi(x, \eta') - (y - x) \varphi_x(x, \eta'),\]

\(\rho\) vanishes at order \(\geq 2\) in \(y = x\). For checking (3.10), we see
\[
Q(a e^{i \frac{\varphi}{h}}) = \int\int e^{i \frac{1}{h} \left((x' - y') \theta' + \varphi(y, \eta')\right)} \frac{q(x, \theta')}{(2\pi h)^{n-1}} a(y, \eta') \, dy' d\theta'. \tag{3.11}
\]
a having compact support and \(\varphi\) being a symbol, one has $|\varphi'_y| \leq C' \quad \text{on} \ \text{supp}(q a)$, so for \(|\theta'| \geq C > 0\), large, the phase of (3.11) $H = (x'-y')\theta' + \varphi(y, \eta')$ satisfies 
\[|H'_{y'}| \geq c(1 + |\theta'|), \quad \text{for } |\theta'| \geq C > 0, \text{ large}\] and if we split the integrand of (3.11) in $q a = \chi(\theta') qa + (1 - \chi(\theta')) q a$, we integrate the second term by parts and obtain 
\[
Q(a e^{i \frac{\varphi}{h}}) = Q'_\varphi({a}) e^{i \frac{\varphi}{h}} + R(a), \tag{3.12}
\]
where $Q_\varphi'$ is a \(\mathcal{G}^s\) symbol of order \(\widetilde{S}_s^{m, k+1}\) having the expansion (3.10) by the stationary phase lemma as \(R(a)\) is an \(\mathcal{O}_s(h^\infty)\) remainder. It is easy to see in view of these arguments that \(a_1 \in \widetilde{S}_s^{m-1, k}\). Moreover, it is to be observed that the above expansion is only a formal Gevrey \(2s-1\) symbol.

The microlocal invertibility of FIO reduces to the PDO case. We refer to \cite{Lascar1988} for a proof of the Gevrey elliptic result in classes \(S_s^{m}\).

For proving Theorem 2, we rewrite (3.8) in the form
\[
(h D_{x_1} + Q(x, h D_x; h)) F u = F(h D_{x_1} + Q') u,
\]
close to \((x_0, \xi_0; x_0, \xi_0)\) for some PDO \(Q'(x, h D_x; h)\) of bi-order \((-1, 0)\) in using a left microlocal inverse of \(F\) close to \((x_0, \xi'_0)\). Indeed, we compute \(F F^*\) and \(F^* F\), and one has writing \(y = (x_1, y')\).

One has, following Eskin \cite{Eskin2011},
\[
F F^*{u}(x, h) = \frac{1}{(2 \pi h)^{n-1}} \iint e^{\frac{i}{h} (\varphi(x, \xi') - \varphi(y, \xi'))} a(x, \xi') \overline{a(y, \xi')} u(x_1, y') \, dy' \, d\xi', 
\]
\(\varphi(x, \xi')\) having been obtained in (3.6). We split the integral above into two terms. The first is a \(h\)-PDO, the second is a smoothing operator. First, we note that the map:
\[
(x, y', \xi') \to (x, y',  \Sigma(x, y', \xi')),
\]
with
\[
\Sigma(x, y', \xi') = \int_0^1 \varphi'_{x'} (x_1, y' + t(x' - y'), \xi') \, dt
\tag{3.13}\]

is a \(\mathcal{G}^s\)-diffeo in a neighbourhood of \((x_0, y_0', \eta_0')\) with \(|x_1| \leq \delta\), \(|x' - y'| \leq \delta\), \(0 < \delta\) small, close to the identity.

Let \((x, y', \eta') \to (x, y', \Sigma^{-1}(x, y', \eta'))\) be an inverse map.

One has obviously $\varphi(x, \xi') - \varphi(y, \xi') = \Sigma(x, y', \xi')(x' - y'),$
and
\[
F F^* u(x, h) = K_1 u(x, h) + K_2 u(x, h),\tag{3.14}
\]
\begin{align}
K_1 u(x, h) &= \frac{1}{(2 \pi h)^{n-1}} \iint e^{\frac{i}{h} (x'- y') \eta'} a(x, \Sigma^{-1}(x, y', \eta')) \notag \\ 
&\quad \times \overline{a(y, \Sigma^{-1}(x, y', \eta'))} \chi\left(\frac{x' - y'}{\delta}\right) |J(x,y', \eta')| u(x,y') \, dy' \, d\eta'. \tag{3.15}
\end{align}

with \(\chi \in \mathcal{G}_0^s(\mathbb{R}^{n-1})\), \(\chi \equiv 1\) for \(|t| \leq \delta/2\), \(\mathrm{supp} \chi \subset |t| \leq \delta\),

\[
K_2 u(x, h) = \frac{1}{(2 \pi h)^{n-1}} \iint e^{\frac{i}{h} (\varphi(x, \xi') - \varphi(y, \xi'))} a(x, \xi') \overline{a(y, \xi')} (1 - \chi)\left(\frac{x' - y'}{\delta}\right)u(y, h) \, dy' \, d\xi'. 
\]

In (3.15),$J = \partial_{\eta'} \Sigma^{-1}$ is a Jacobian.

$K_1$ is a \(\mathcal{G}^s\) PDO of order \((2m, 2k)\). Its double symbol may be reduced to a single one, its principal symbol is $|a_0(x, \Sigma^{-1}(x, y, \eta'))|^2 \, |\varphi_{x', \xi'}''(x, \Sigma^{-1}(x, y, \eta'))|^{-1} $ with \(a_0\) the principal symbol of \(F\), \(|\varphi_{x', \xi'}''|\) is the absolute value of \(\det \varphi_{x', \xi'}''\). \(K_2\) is a \(\mathcal{G}^s\)-smoothing operator. Its kernel is in \(\mathcal{G}^s\) and has a decay $e^{\frac{-1}{C}h^{-1/s}}$ in $\mathcal{G}^s_b(\mathbb{R}^{n}\times \mathbb{R}^{n-1})$.

Let us outline it. On \(\mathrm{supp} (1 - \chi)\),  $|(\varphi'_{\eta'}(x, \xi') - \varphi'_{\eta'}(y, \xi'))| \geq c > 0 $ and the non-stationary phase lemma applies. We perform IPP using $L = \frac{h}{i} \, (\varphi'_{\eta'}(x, \xi') - \varphi'_{\eta'}(y, \xi')) / |\cdot|^2 \partial_{\xi'}$ and Gevrey quasi-norms (see Lemma 2.1).

Let \( F \) be a Fourier Integral Operator, one has to compute
\[
F^*u(x,h) = (2\pi h)^{-(n-1)} \int \overline{a}(y, \xi') e^{i/h (x' \xi'- \varphi(x_1, y', \xi'))} u(x_1, \eta') \, d\xi' \, dy',
\]
and one has
\begin{align}
F^*Fu(y,h) &= \frac{1}{(2\pi h)^{2(n-1)}} \int \overline{a}(x, \xi')a(x, \eta') e^{\frac{i}{h} (y' \xi'- \varphi(y_1, x', \xi'))} \notag\\
&\quad \times \varphi(y_1,x',\eta')\hat{u}(y_1, \frac{\eta'}{h'}) \, dx' d\xi' \, d\eta'. \tag{3.16}
\end{align}

We reduce the integration in \((3.16)\) to variables \(d\eta'\) by  stationary phase. 
Indeed, let us set
\[
I_a(y, \eta') = \frac{1}{(2\pi h)^{n-1}} \int \overline{a}(x, \xi') a(y_1,x', \eta') e^{\frac{i}{h} \psi_{x, \eta''}} dx' d\xi', \tag{3.17}
\]
then \(I_a\) is:

\[
I_a = A(y, \eta') e^{\frac{i}{h} \psi_{y', \eta'}(y_0', \eta')} + R_a, \quad \text{with } R_a \in \mathcal{O}_s(h^\infty) \text{ and with some }
\]
\[
A \in \widetilde{S}_s^{2 m', 2k'} \quad \text{if} \quad  a \in \widetilde{S}_s^{m',k'} \quad \text{with main part } 
\]
\[
|a_0(y_1, y_0',\eta')|^2 |\varphi_{x', \xi'}''(y_1, y_0', \eta')|^{-1},
\]
where \( y_0' \) satisfies \( \varphi_{\eta'}'(y_1, y_0',\eta') = y' \).

For obtaining \((3.17)\), we split \( I_a \) into two parts, writing
\[
1_{\xi'} = \chi\left(\frac{\xi'-\eta'}{|\eta'|}\right) + (1 - \chi)\left(\frac{\xi'-\eta'}{|\eta'|}\right),
\]
with \( \chi \in \mathcal{G}_0^s(\mathbb{R}^{n-1}) \), \( \text{with support on }|t| \leq 1/2,  \chi \equiv 1 \text{ for }|t| \leq \frac{1}{4}\}, \) and
\[
\psi_{y, \eta'}(x', \xi') = y' \xi' - \varphi(y_1, x', \xi') + \varphi(y_1, x', \eta').
\]
Satisfies
\[
\partial_{x'} \psi_{y, \eta'} = 0 \implies \varphi_{x'}'(y_1, x', \eta') - \varphi_{x'}'(y_1, x', \xi') = 0,
\]
which implies \(\xi' = \eta'\), since
\[
|\partial_{x'}\psi_{y,\eta'}|\geq c |\xi'-\eta'|, \quad \text{since } \det \partial_{x', \eta'}^2 \varphi \neq 0.
\]
However, if
\[
H = 
\begin{pmatrix}
0 & -\varphi_{x', \xi'}''(x,\xi') \\
-\varphi_{x', \xi'}''(x,\xi') & -\varphi_{\xi', \xi'}''(x', \xi')
\end{pmatrix}
\]
one has
\[
|\det H| = |\det \varphi_{x', \xi'}^{''}|^2, \quad \text{sgn } H = 0,
\]
and so \((\psi_{y, \eta'})_c = y' \eta';\) the stationary phase applies, and one has
\[
A \in \widetilde{S}_s^{ 2m', 2k'}, \quad \text{since } |\xi'| \sim |\eta'| \text{ on supp } \chi\left(\frac{\xi'-\eta'}{|\eta'|}\right).
\]
Now,
\[
R_a(y, \eta') = \int \overline{a}(x, \xi') a(y_1, x',\eta') e^{(i/h) \psi_{y, \eta', \xi'} (x', \eta')} (1 - \chi) \left(\frac{\xi'-\eta'}{|\eta'|}\right) \frac{1}{(2 \pi h)^{n-1}}d\xi' dx',
\]
for obtaining \(R_a \in \mathcal{O}_s(h^\infty)\), we perform IPP with
\[
L = \frac{h}{i} \frac{\partial_{x'} \psi_{y, \eta'}}{|\partial_{x'} \psi_{y, \eta'}|^2},
\frac{\partial}{\partial x'}\]
since \(|\partial_{x'} \psi_{y, \eta'}| \geq c |\xi' - \eta'| \geq c' |\eta'| \) on \(\text{supp}(1-\chi)\).

For details, we refer to Lemma 2.1 of these notes and to \cite{LascarLascar1997}.

Gathering the above remarks we have proved that \( F^* F\) is a \( G^s \) PDO of order \( (2m',2k') \) if \( a(x,\eta') \in \tilde{S}_{s'}^{m',k'} \).

Now for obtaining (3.11), it suffices to construct \( G^s \) microlocal parametrices to elliptic PDOs at \( (x_0, \xi_0) \) given. This is done in \cite{Lascar1988}.

At the next step, we construct \( G^s \) PDOs \( A_2, B_2, \text{s.t.}\) 
\[
B_2 (hD_{1} + Q') A_2 - hD_{x_2}, \quad B_2 A_2 - \text{Id}, \quad A_2 B_2 - \text{Id are small } \mathcal{G}^s \text{ remainders} \tag{3.18}.
\]

\underline{Step 2}

For achieving (3.18), we use the PDO calculus sketched above and an induction using quasi-norms.

We proved above if \( q \in \tilde{S}_s^{m_0,k_0} \), \( a \in \tilde{S}_s^{m',k'} \) one has 
\[
(q o a)(x,\xi') = (q a)(x,\xi') + r_1(q,a)(x,\xi') 
\tag{3.19}\]
with 
\[
r_1(q,a) \in \tilde{S}_s^{m_0+m'_{1},k_0+k'-1}
\]

The following lemma is direct.

\subsection*{{Lemma 3.1}} Let \(q \in \tilde{S}_s^{{-1},{0}}\) (\(\mathbb{R}^n \times \mathbb{R}^{n-1}\)) i.e.\ \(q\) is of \(bi\)-order (-1,0). There is a sequence \(\bigl(a_j\bigr)_{j\ge0}\) of Gevrey symbols of bi-order \(\bigl({-j},{-j}\bigr)\), i.e.\ \(a_j \in \tilde{S}_s^{-j,-j}\) such that

\[
\begin{cases}
\tfrac{1}{i}\,\partial_{x_1}a_0 \;+\;\tfrac{1}{h}\,q\,a_0 \;=\; 0,\\[6pt]
a_0\bigl\lvert_{x_1=0} = 0 \quad \text{elliptic at }(x_0,\eta_0'),
\end{cases}
\quad
\begin{cases}
\tfrac{1}{i}\partial x_1\,a_j \;+\;q\,\tfrac{ a_j}{ h} \;=\;-\tfrac{r_1\bigl(q,a_{j-1}\bigr)}{a},\quad j \ge 1\\
a_j|_{x_1=0}=0
\tag{3.20}
\end{cases}
\]

\textbf{Proof.} Equations (3.20) are solved by induction,

\[
a_0 \;=\;\exp\Bigl(-\tfrac{i}{h}\int_{0}^{x_1}qds\Bigr), 
\quad
a_j \;=\; -\,i\,\frac{a_0}{h} \;\int_{0}^{x_1}\,\tfrac{r_1\bigl(q,a_{j-1}\bigr)}{a_0}\,ds
\quad (j\ge1).
\]

So one has \(a_j \in \tilde{S}^{-j,-j}_s\) and if \(A_j = \mathrm{op}_h\bigl(a_j\bigr)\) one has for \(N>0\) if

\[
A^{(N)} = A_0 + \cdots +A_{N-1}, \quad (h D_{x_1} + Q) A^{(N)} = h A^{(N)} D_{x_1} + R_N, \quad R_N \in \text{PDO of class } \tilde{S}_s^{-N-1, -N}.
\]

\(A^{(N)}\) being elliptic at \((x_0, \eta_0')\), we have reduced the problem to proving that \(h D_{x_1} + Q\), with $d^{\circ}  Q = (m_0 -1, m_0), m_0 \leq -n$ is conjugate to \(h D_{x_1}\) microlocally. Since \(h D_{x_1} + Q\) is conjugate to \(P\), we will have \(h D_{x_1}\) and \(P\) conjugate at \((x_0, \xi_0)\)

In (3.19), we split \(r_1(q, a)\) into two terms, writing in the integral giving
\[
1 = \chi\left(y',\frac{\eta'}{|\xi'| }\right) + (1 - \chi)(y',\frac{\eta'}{|\xi'| }),
\]
with \(\chi \in \mathcal{G}^s_0(\mathbb{R}^{2n-2})\), \(\chi\) having its support close to \(0\) and being \(1\) close to \(0\).

We write
\[
r_1(q, a) = R_\chi(a) + S_\chi(a), \tag{3.21}
\]
where \(S_\chi\) is a Gevrey small remainder by non stationary phase lemma (see \cite{Lascar1988} and \cite{Hitrik2023complex}).

Assume \(q\) is of bi-order \((m_0 -1, m_0)\) with \(m_0 \leq -n\), given. Define \([a_j]\) as a sequence by
\begin{equation}
\left\{
\begin{aligned}
    &\frac{h}{i} \partial_{x_1} a_0 + q a_0 = 0, \\
    &a_0 |_{x_1 = 0} \text{ elliptic, given in } \tilde{S}_s^{0,0}
\end{aligned}
\right.
\quad
\left\{
\begin{aligned}
    &\frac{h}{i} \partial_{x_1} a_j + q a_j = -R_\chi(a_{j-1}), \\
    &a_j |_{x_1 = 0} = 0 \text{ for } j \geq 1.
\end{aligned}
\right.
\tag{3.22}
\end{equation}

we have:

\subsection*{{Lemma 3.2}} If \(m_0 \leq -n\), the above (3.22) sequence is a formal Gevrey \(s\) symbol of bi-order \(0\), i.e., one has for some \(C > 0\),
\[
|\partial_x^\alpha \partial_{\xi'}^\beta a_j(x, \xi')| \leq C^{ (1  + |\alpha| + |\beta|+j)} \, \alpha!^s \, \beta!^{s} \,j^{s}! \, h^j\langle \xi \rangle^{-|\beta|-j}\tag{3.23}
\]
for all \(j \geq 0\), \(\alpha \in \mathbb{N}^n\), \(\beta \in \mathbb{N}^{n-1}\), \((x, \xi') \in \mathbb{R}^n \times \mathbb{R}^{n-1}\), and \(h \in (0, 1]\). 

\textbf{Proof.} The PDO calculus guarantees \( R_\chi(a) \in \tilde{S}^{m_0 - 1+ m',k'-1+m_0} \), if \( q \in \tilde{S}^{m_0, m_0} \), \( a \in \tilde{S}_s^{m', k'} \). In order to prove \((3.23)\), we will proceed by induction using Gevrey quasi norms. Let us set:
\[
\tilde{N}_{m', k'}(a, T) = \sum_{(\alpha, \beta)} \frac{|\partial_x^\alpha \partial_\xi^\beta a(x, \xi')|}{\alpha!^s \, \beta^s!} \frac{T^{|\alpha|+|\beta|}}{\langle \xi' \rangle^{k'-|\beta|} h^{-m'}}, \tag{3.24}
\]
and
\[
\overline{N}_{m', k'}(a, T) = \sup_{(x, \xi') \in \mathbb{R}^n \times \mathbb{R}^{n-1}} \tilde{N}_{m', k'}(a, T)(x, \xi').
\]
We note that by the Leibniz formula:

\[
\overline{N}_{m+m', k+k'}(aa', T) \leq \overline{N}_{m, k}(a, T) \, \overline{N}_{m', k'}(a', T). \tag{3.25}
\]

One may prove also that for \(\epsilon > 0\) small, \(T > 0\) small, and all \((\alpha, \beta) \in \mathbb{N}^n \times \mathbb{N}^{n-1}\), the bound
\[
\overline{N}_{m,k -|\beta|}( \partial_x^\alpha \partial_\xi^\beta a, T) \leq \epsilon^{-s(|\alpha| + |\beta|)} T^{(-|\alpha| - |\beta|)} \, \alpha!^s \, \beta!^s \, \overline{N}_{m, k}(a, T(1 + \epsilon)^s). \tag{3.26}
\]

Now, if \(q \in \tilde{S}^{m_0, m_0}\), \(m_0 \leq -n\), one has quite obviously from \((3.25)\), \((3.26)\), the estimates for some $M_0(q)>0$:
\[
\overline{N}_{m'-1, k'-1}(\frac{1}{h} R_\chi(a), T) \leq C_n r^{2(n-1)} \ \epsilon^{-s} {T^{-1}} \, M_0(q)\overline{N}_{m', k'}(a, T(1 + \epsilon)^{s}) \tag{3.27}
\]

If \(a \in \tilde{S}_s^{m'; k'}\) and if \(R_\chi(a)\) is defined by \((3.21)\),\((3.27)\) is rough, but the precise control of quasinorms of \({J^h}(b)\) with respect to the constants of \(b\) is not easy. So we prove by induction on \(j \geq 0\) that for some \(M'_0(n, q) > 0\), \(\epsilon > 0\) small, \(T > 0\) small, the bounds
\[
\overline{N}_{-j, -j}(a_j, T) \leq {M'^{j}_0 {T}^{-j}  \epsilon^{-s j} j^{s j} \overline{N}_{0}(a_0, 2T)} \tag{3.28}
\]
hold. 

Indeed, we check that for \(j \geq 1\), for all small \(\epsilon > 0, T > 0\), one has the bounds
\[
\overline{N}_{-j, -j}(a_j, T) \leq {M'^{j}_0 {T}^{-j}  \epsilon^{-s j} j^{s j} \overline{N}_{0}(a_0, T(1+\frac{\epsilon}{j})^{sj})} \tag{3.29}
\]

From (3.27), (3.29) is true for \(j = 1\), assuming (3.29) for \(j - 1, j \geq 2\) for \(\epsilon > 0\), $\tilde{\epsilon} > 0$ small, one has:
\[
\overline{N}_{-j, -j}(a_j, T) \leq \epsilon^{-sj} M'^j_0 \, T^{- j} (1 + \epsilon)^{-s(j - 1)} (j-1)^{s(j  - 1)} \overline{N}_0(a_0, (1 + \epsilon)^s({(1 + \frac{\tilde{\epsilon}} {j-1})^{s( j - 1)}}T). \tag{3.30}
\]

Fixing \(\delta > 0\) small, we choose in (3.30): $\epsilon = \delta / j, \quad \frac{\tilde{\epsilon}}{{j-1}} = \delta / j$, so one has
\[
(1 + \epsilon)^s (1 + \frac{\tilde{\epsilon}}{j-1})^{s(j-1)} \leq (1 + \delta / j)^{sj}.
\]

And for proving (3.28), one has to check:
\[
\left( 1 + \epsilon \right)^{-s( j-1)} \tilde{\epsilon}^{-s j} \delta^{s j} {{(s-1)^{s(j-1)}}} \leq \left(\frac{\delta}{j}\right)^s. \tag{3.31}
\]

(3.31) can be rewritten as:
\[
\left(1 + \frac{\delta}{j}\right)^{-s( j-1)} \left(\frac{\delta(j-1)}{j}\right)^{-s (j-1)} \left(j-1\right)^{s( j-1)} \left(\frac{\delta}{j}\right)^{sj}\leq \left(\frac{\delta}{j}\right)^s.
\]

Which is:$\left(1 + \frac{\delta}{j}\right)^{-s( j-1)} \leq 1,$  which is true since \(j \geq 1\). So (3.29) is proven. From it one has:\[
\overline{N}_{-j, -j}(a_j, T) \leq M'^j_0 \, T^{-j} \delta^{-s j} j^{s j} \overline{N}_0(a_0, 2T), \tag{3.32}
\]
since 
\[
\log(1 + \delta / j)^{sj} \leq \delta j \quad \text{and} \quad (1 + \delta / j)^{sj} \leq 2 \quad \text{for} \quad \delta > 0 \, \text{ small}.
\]

(3.32) means precisely that \((a_j)_{j \geq 0}\) is a formal Gevrey symbol in \(\widetilde{S}_s^{0,0}\).

This ends Lemma 3.2's proof.

Now, by the Carleson moment theorem \cite{Carleson1964}, adapted by Boutet de Monvel-Kreé \cite{BoutetDeMonvel1967}, we obtain a realization of the sequence \((a_j)_{j \geq 0}\), \(a \in \widetilde{S}_s^{0,0}(\mathbb{R}^n \times \mathbb{R}^{n-1})\), with 
\[
\left| \partial_x^\alpha \partial_{\xi'}^\beta (a - \sum_{j < N} a_j) \right| 
\leq C ^{1+|\alpha| + |\beta|+N}\, N!^{s}\alpha!^{s} \beta!^sh^N \langle \xi' \rangle^{-|\beta| - N}, \tag{3.33}
\]

Setting \(A = \operatorname{op}_h a\), \(A\) solves \([hD_{x_1},A]+QA\) is a \(\mathcal{G}^s\) remainder.

Indeed, it has symbol $\frac{h}{i} \partial_{x_1} a + q a + r_1 (q, a) = R_\chi (a_{N-1}) + \left(\frac{h}{i} \partial_{x_1} + q + R_\chi \right)\rho^N+S_x(a),$
with \(\rho^N = a - \sum_{j < N} a_j\), and in view of (3.33), one has the result. This ends the proof of Theorem 2. Let us deduce now Theorem 1.
\section{Proof of Theorem 1}
We start by considering a zero-order PDO reduced by the implicit function theorem to the form:
\begin{equation}
P = hD_{x_1} + Q(x, hD_{x'}; h),
\tag{4.1}
\end{equation}
where \( Q\) is a Gevrey PDO on \(\mathbb{R}^{n-1}\) of zero bi-order, depending on \(x_1\) as a parameter having split the variable \(x\) of  \(\mathbb{R}^n\) into: $x = (x_1, x') \quad \text{with } x_1 \in \mathbb{R}, \; x' \in \mathbb{R}^{n-1},$
and may also assume that:
\begin{equation}
S = \{x_1 = 0\}.
\tag{4.2}
\end{equation}

Let $P$ be given by (4.1) using FIOs we reduce $P$ to $P'$

$P' = h D_{x_1} + Q'(x,h D_{x'},h)$ with $Q'$ a $G^s$ PDO of bi-order -1 ( i.e $\sigma_{Q'}(x,\xi') \in \Sigma_s^{-1,-1} ( \mathbb{R}^n \times \mathbb{R}^{n-1} )$. Using Kuranishi's trick we reduce the WKB problem (1.6) to the same question with $\varphi \equiv 0$. Let us note that in that case there is no need of almost holomorphic extension of the symbols since the function $\varphi$ at the beginning is assumed real ( see (1.5) ) and in that case Kuranishi's trick is a reduction of a double symbol of a PDO to a single one.

If we solve (4.3) below:\[
\begin{cases}
\dfrac{h}{i}\,\partial_{x_1}\,a \;+\; e^{-\tfrac{i\,\phi_0}{h}}\,
Q'\bigl(x,h\,D_{x'},h\bigr)\Bigl(e^{\tfrac{i\,\phi_0}{h}}\,a\Bigr)
\;=\; 0,\\[6pt]
a\bigl|_{x_1=0} \text{ is elliptic, with } \phi_0(x,\xi') = x'\,\xi',
\end{cases}
\tag{4.3}
\]

We will solve the homogeneous WKB problem by setting \( a(x; h) = a(x, 0; h) \), with \( a(x, \xi'; h) \) as a solution of (4.3) 

Solving the inhomogeneous WKB problem is then standard. Moreover,  in (4.3) \( \equiv 0  \) means modulo small \( \mathcal{G}^s \)-remainders \( \mathcal{O}_s(h^\infty) \).

Now solving (4.3) is equivalent to finding \( A(x, hD_{x'}) = \text{op}_h(A) \) elliptic, such that:
\begin{equation}
P'(x, hD_{x}) A \equiv h A D_{x_1}, \quad \text{close} (x_0, \xi_0).
\tag{4.4}
\end{equation}

Equation (4.4) is a consequence of Theorem 2, so the proof of Theorem 1 is finished. A corollary of Theorem 1 and Theorem 2 could be the following:

\paragraph{Corollary 3.5.}
Let \( P(x, hD_x; h) \) be an evolution equation:
\[
P = h D_{x_1} + Q(x, h D_x'),
\]
where \( Q \) is a zero-biorder PDO on \( \mathbb{R}^{n-1} \), with a real principal type close to \( (x_0, \xi_0) \).

There exists an \( h \)-FIO of the form:
\[
A u(x, h) = (2 \pi h)^{-(n-1)} \int e^{i S(x, \eta') / h} a(x, \eta') \hat{u}(x_1,\eta'/ h) \, d\eta,
\]
such that the family \( u(x, h) = A(u_{0} \otimes 1_{x_1}) \) is a singular Gevrey s solution for \( P \). That is, the Gevrey wavefront set \( WF_{s,h}(u(x, h)) \) consists of a segment of the null bi-characteristic of \( p(x, \xi) \) through \( (x_0, \xi_0) \). \( u \) is obtained in taking for \( u_0 \) a coherent state at \( (x_0', \xi'_0) \in \mathbb{R}^{n-1}_{x'} \times \mathbb{R}^{n-d}_{\xi'} \), i.e.,
\begin{equation}
u_0(x', h) = (\pi h)^{-(n-1)/4} e^{\xi'_{0} \cdot i(x' - x'_0) / h} e^{- \frac{|x - x'_0|^2}{2h}}.
\tag{4.5}
\end{equation}

One may easily check \( WF_h(u_0) = WF_{s,h} (u_0) = { (x_0', \xi'_0)}, \) for all \( s \geq 1\)  using the Bargman transform:

\begin{equation}
T_0 u(x',h) = C_n h^{-3(n-1)/4} \int_{\mathbb{R}^{n-1}} e^{\frac{i}{h} \varphi_0( x', y')}u(y', h) dy'
\tag{4.6}
\end{equation}

From \cite{Zworski2012}, we see

\(
T_0 u(x',h) = C_n h^{-3(n-1)/4} e^{-\phi_0 (x')/h}  e^{-Re (x') Im(x')/h}\mathcal{F_h}(ue^{(\cdot -Re(x'))^2/2h})(-Im(x')))\)

So computing \( (T_0 u_0)(x', h) \), one may prove

\begin{equation}
\left| T_0 u_0(x', h) \right| = e^{\frac{1}{h}\phi_0(x')} C'_n h^{-(n-1)} e^{-\frac{1}{2h} |x' - x'_0|^2} \text{   for  some } C'_n >0 
\tag{4.7}
\end{equation}

Hence, \( WF_{s,h}(u_0) \) is computed.

Note that improvements could be obtained in looking at images of Gaussian beams by real FIO’s in the spirit of Cordoba-Fefferman works.

\section{Proof of Theorem 3}

\textbf{Step 1.} Construction of the canonical transform.

The construction  of \( \kappa : T^*\mathbb{R}^n \setminus0 \to \Lambda_{\phi_0} \) transforming \( p(y,\eta) \) into \( \xi_1|_{\Lambda_{\phi_0}} \) is easy. At first we construct by Darboux lemma $\chi$ a real $\mathcal{G}^s$ canonical transformation near \( (y_0, \eta_0) \), mapping \( p(y,\eta) \) to \( \eta_1 \), here \( p(y,\eta) \) is the principal symbol of \( P(y, D_y) \). Then we compose \( \chi \) by \( \kappa_{T_0}| _{\mathbb{R}^{2n}} \) and set

\[
\kappa =  \kappa_{T_0}|_{\mathbb{R}^{2n}} \circ \chi \tag{5.1}
\]

The points \( (x, \xi) \in \Lambda_{\phi_0} \) are parameterized by \( x \) since \( \xi(x) = \frac{2}{i} \partial_x \phi_0(x) \), i.e.  $\xi(x) = -Im(x)$ We set
\[
x_0 = (y^0_1,y'_0-i\eta'_0) \in \mathbb{C}^n, \quad \text{close to $x_0$ we set} \quad (y(x), \eta(x)) = \kappa^{-1}(x, \frac{2}{i} \partial_x \phi_0(x')).
\]
Then the map \( x \in \mathbb{C}^n \to y(x) \in \mathbb{R}^n \) is of class \( \mathcal{G}^s \) and \( dy \) the total differential of \( y \) is surjective $T_x \mathbb{C}^n \to T_{y(x)} \mathbb{R}^n$, since \( \text{Re} \left( \frac{\partial y}{\partial \bar{x}} \right) \) is bijective as \( d\kappa \) is a symplectic linear transform (a proof of this is given below). So we may define
\[
\Gamma = \{ ( x, y(x) ) \in \mathbb{C}^n \times \mathbb{R}^n : x \in \text{neigh}(x_0, \mathbb{C}^n) \}. \tag{5.2}
\]
\(\Gamma\) is then a submanifold of \( \mathbb{C}^n \times \mathbb{C}^n \), and \( \Gamma \) is totally real i.e. \( T\Gamma \cap iT\Gamma = \{o\}\), and $\Gamma$ is of maximal dimension, i.e. \( dim_{\mathbb{R}}\Gamma = 2n \).

Using local coordinates near \( \Gamma \), we may construct for smooth functions on \( \Gamma \) extensions to \( \mathbb{C}^n \times \mathbb{C}^n \) which are almost holomorphic w.r.t. \( \Gamma \) i.e. with their $\overline{\partial}$-flat on \( \Gamma \). It applies in particular to \( \mathcal{G}^s \) maps and for such maps \( f \) there are extensions \( \tilde{f} \) to \( \mathbb{C}^n \times \mathbb{C}^n \) in the same \( \mathcal{G}^s \)-class and with $\left( \partial_{\overline{x}} \tilde{f}, \partial_{\overline{y}} \tilde{f} \right) = \mathcal{O}\left( \exp \left( -\frac{1}{C} d(\cdot, \Gamma)^{\frac{-1}{s-1}} \right) \right)$  (see \cite{Hitrik2023complex} for a proof).

\textbf{Step 2.} \underline{Eikonal Equation}

The statement \( (1.17) \) is equivalent to solving eikonal and transport equations for the determination of the phase and the symbol \( a \) of the FBI map $T$.
The eikonal equation is
\[
\varphi'_{x_1} = p(y,-\varphi'_y) \quad \text{on } \Gamma 
\tag{5.3}
\]
It means we want to construct an almost holomorphic map \( \varphi \) w.r.t. \( \Gamma \) near \( (x_0, y_0) \in \mathbb{C}^n \times \mathbb{C}^n \) such that the holomorphic derivatives of \( \varphi \), \( (\varphi'_x, \varphi'_y) \) are real on \( \Gamma \) and satisfy (5.3) on \( \Gamma \).
Let us note \( \psi \) the map defined on \( \Gamma  \approx  \mathbb{C}^n \) obtained in solving
\[
d\psi = \xi(x) dx - \eta(x) dy|_\Gamma
\]

The 1-form  $ \xi dx - \eta dy|_\Gamma$ is closed near  $x_0$ since $dx$ is canonical. $\psi$  is well defined up to an additive constant. It is a $ \mathcal{G}^s$  map since  $\kappa$  is a one . We note $\varphi$ an almost holomorphic extension w.r.t. $\Gamma$. One has:

\[
\varphi(x, y(x)) = \psi(x), \quad \partial_{\overline{x}} \varphi = \partial_{\overline{y}} \varphi = 0 \quad \text{on } \Gamma.
\]
\[
d\psi(x) = \left(\xi(x) - \text{  }^{t}(\frac{\partial y}{\partial x}) \right) \eta(x) dx - \text{  }^t(\frac{\partial y}{\partial \overline{x}}) \eta(x) d \overline{x}  \tag{5.4}
\]
\[
\text{Since Re} {\text{ }^t}(\frac{\partial y}{\partial x}) \text{ is bijective, one has on } \Gamma
\]
\[
\frac{\partial \varphi}{\partial x} = \xi(x), \quad \frac{\partial \varphi}{\partial y} = - \eta(x), \quad x \in \text{neigh}(x_0, \mathbb{C}^n) \tag{5.5}
\]
So (5.3) is satisfied near \( (x_0, y_0) \) since
\[
\varphi_{x_1} = \xi_1(x)= - Imx_1 =  p(y(x), \eta(x)) = p(y(x), -\varphi'_{y}(x,y(x)) \tag{5.6}\]
\textbf{Step 3.} Transport equation.  
The phase function \( \varphi(x, y) \) is constructed near \( (x_0, y_0) \). Let us begin to construct an amplitude \( a(x, y; h) \) solving

\[
h D_{Rex_1} T_{a}u - h^m T_{a} P u \sim 0 \quad \text{in} \quad \mathcal{G}^{2s-1} \text{ close to } x_0. \tag{5.6}
\]

We need to solve

\[
P'a = e^{-i \varphi / h} \left\{ h D_{Rex_1} - h^{m t}{P}(y, D_y) \right\} \left( a e^{i \frac{\varphi}{h}} \right) \equiv 0 \tag{5.7}
\]

is \( \mathcal{O}_s(h^\infty) \) for a symbol \( a(x,y,h) \), elliptic at \((x_0, y_0)\) of \( \mathcal{G}^s \) class which will be chosen as the amplitude of \( T\).

First we note that \( h^{mt}P(y, D_y) \) is an \( h \)-PDO with principal symbol  $\tilde{p}(y, \eta) = p(y, -\eta)$, it belongs to $\tilde{S}^{0,m}_s$. We set
\begin{equation}
P_{\varphi}(x, y, hD_y) = e^{-i \varphi/ht} P(y, D_y) \,h^m \left( a e^{i \varphi/h} \right) 
\tag{5.8}
\end{equation}

$P_{\varphi}$ is an h PDO of degree  (0, m), its principal part is $p(y, -\eta - \varphi'_y)$. So \( P' \) defined in \( (5.7) \) is an \( h \) PDO of biorder \( (0, n) \) with principal symbol 
\begin{equation}
p'(x, y, \xi, \eta) = \text{Re} \xi_1 + \varphi' _{\text{Re} x_1} - P_\varphi(x,y,\eta) .
\tag{5.9}
\end{equation}
One has to solve a WKB problem for \( P' \) with \( \varphi' \equiv 0 \). So the eikonal equation is satisfied on $\Gamma$ since
$P'|_{\xi = \eta = 0} = \varphi'_{x_1} - p(y, - \varphi'_y) \, \text{on} \, \Gamma \quad \text{which is 0 by (5.3)}.$ Note that since \( \varphi \) is complex valued, we have taken in \( (5.8) \) an almost holomorphic extension of \( p(y,\eta) \) to \( \mathbb{C}^n \times \mathbb{C}^n \), which is done obviously for differential operators in taking extensions of the coefficients.

We set now:
\begin{equation}
\alpha' = \frac{\partial p'}{\partial \xi} \Big|_{\xi = \eta = 0}, \quad \beta' = \frac{\partial p'}{\partial \eta} \Big|_{\xi = \eta = 0}
\tag{5.10}
\end{equation}

One has \( \alpha' = e_1 \), \( \beta' = \frac{\partial p}{\partial \eta} \left( y_1, - \varphi'_y \right) \) on \( \Gamma \).

Consider:
\begin{equation}
V' = a' \partial_x + \overline{a'} \partial_{\bar{x}} + \bar{b'} {\partial \bar{y}} +{b' \partial y}
\tag{5.11}
\end{equation}
a real vector field on $\mathbb{C}^n \times \mathbb{C}^n$. One checks \( V'|_\Gamma \) is tangent to \( \Gamma \) iff,
\begin{equation}
b' \in \mathbb{R}^n, \quad b' = \frac{\partial y}{\partial x} \, a' + \frac{\partial y}{\partial {\bar{x}} } \, \bar{a'} \quad \text{on} \quad \Gamma
\tag{5.12}
\end{equation}

If \( V = \alpha {\partial x} + \bar{\alpha} {\partial \bar{x}} \) is a vector field in \( \mathbb{C}^n \simeq \Gamma \), one has $V'u' = V u \quad \text{on} \quad \Gamma$
and if \( \alpha = a' \), $b'=\frac{\partial y}{\partial x}\alpha' +\frac{\partial y}{\partial \bar{x}}\bar{\alpha}'$ and if $\partial_{\bar{x}}  u' = \partial_{\bar{y}} u' = 0$ on $\Gamma$.
With
$b' = \frac{\partial p}{\partial \eta} (y_1 - \varphi'_y), \quad a' = \frac{\partial p'}{\partial \xi} \Big|_{\xi = \eta = 0}$
(5.12) is satisfied on $\Gamma$. Indeed by design $\chi' H_p = H_{\eta_1} = (e_1, 0)$, setting $\chi'^{-1} = \begin{bmatrix} P & Q \\ R & S \end{bmatrix}$, where $\begin{bmatrix} P & Q \\ R & S \end{bmatrix} $ is real symplectic matrix with P invertible, $^tPR, ^{t}QS$ symmetric $^tPS- ^tRQ$= I. One has:
\[\frac{\partial y}{\partial x} = \frac{1}{2} (P+ iQ),\quad \frac{\partial y}{\partial \bar{x}} = \frac{1}{2} (P- iQ),
\tag{5.13}
\]
\[
\left( \frac{\partial y}{\partial x} + \frac{\partial y}{\partial \bar{x}} \right) e_1 = P e_1 = \frac{\partial p}{\partial \eta} (y, - \varphi'_y) \quad \text{on} \quad \Gamma \text{ since  } \chi^{-1}(Rex, -Imx) = (y(x),\eta(x)),
\]
hence equation (5.12) is satisfied on \( \Gamma\).
We will solve transport equations on \( \Gamma \) and extend solutions almost holomorphically wrt \( \Gamma \). At first, we construct a formal solution \( (a'_j)_{j \geq 0}, a'_j \in S^{-j} \left( \mathbb{C}^n \times \mathbb{C}^n \right) \) and realize it by Carleson's method.

In order to justify transport equations, it seems good to prove some expansions in the complex case.

It is known if \( Q (y, hD_y) \) is a zero order \( h \) PDO on \( \mathbb{R}^n \), \( \varphi \) a complex valued function with
\[
\text{Im} \varphi \geq 0, \quad d\varphi \neq 0 \quad \text{when} \quad \text{Im} \varphi = 0, \quad
\tag{5.14}
\]
 a $a$ symbol on $\mathbb{R}^n $ of degree 0. The expansion
\[
Q \left( a e^{i \varphi} \right) \simeq e^{i \varphi/h} \sum_{\alpha \in \mathbb{N}^n} \frac{h^{|\alpha|}}{\alpha!} \tilde{q}^{(\alpha)} \left( y, \varphi'_y \right) D^{\alpha}_z \left( e^{\frac{i\rho}{h}} a \right) \Big|_{z = y}
\tag{5.15}
\]
with
\[
\rho(y, z) = \varphi (z) - \varphi (y) - (z- y) \varphi'_y (y)
\]
In equation (5.15), we have defined by \( \widetilde{q} \) an almost holomorphic extension of \( q(y, \eta) \text{ to }  \mathbb{C}^n \times \mathbb{C}^n \) is obtained by Melin and Sjöstrand in \cite{HelmSjoberg1975}. We adapt their arguments and compute the remainders. We deal with phase
\[
a(y, X) = (y - z)\eta + \varphi (z)
\tag{5.16}
\]
Here, \( X = (z, \eta) \in \mathbb{R}^{2n}, \varphi \) is a $\mathcal{G}^{s}$ function defined close to \( y_0 \in \mathbb{R}^n, \varphi \) is complex valued with \( \text{Im} \, \varphi \geq 0 \), and \( d \varphi (y_0) = -\eta_0 \in \mathbb{R}^n \setminus 0 \). We will work close to $X_0 = (y_0, \varphi'_y(y_0))$  and we have 

\[{Im} \, a(y, X) \geq 0, \, a'_X(y_0, X_0) = 0, \, \text{and} \,\text{det}  a''_{X,X}(y_0, X_0) \neq 0 \tag{5.17}.
\]

We have the lemma.

\subsection*{{Lemma 5.1}} Let \( u_h(X) \) be a smooth $\mathcal{G}^{s}_0$ symbol of degree 0 in \( \ \mathbb{R}^{2n} \) with compact support close to \( X_0 \), let \( v_h(y) \) be

\[
v_h(y) = e^{-i  a(y, Z_y)/h} \int_{\mathbb{R}^{2n}} e^{\frac{i}{h}a(y, X)} u_h(X) \, dX\tag{5.18}.
\]

\( v_h \) is defined close \( y_0 \in \mathbb{R}^n \) if we have chosen an almost analytic extension of \( a \) to \( \mathbb{C}^{2n} \) of $\mathcal{G}^s$ class and defined \( Z(y) \) as the solution of \( a_Z'(y, Z) = 0 \), \( Z(y_0) = X_0 \), which is possible since

\[
\det a_{Z,Z}''(y_0, X_0) \neq 0 \quad \text{and} \quad a_{Z,\bar{Z}}''(y_0, X_0) = 0.
\]

\( v_h \) is a $\mathcal{G}^s$ symbol of degree \( -n \), it has an expansion

\[
v_h(y) \sim \sum_{\nu=0}^{\infty} h^{\nu+n} C_{\nu,y}(D) u_n(Z_y) \tag{5.19}
\]

with \( C_{\nu,y} \) differential operators of degree \( \leq 2\nu \) and $u_h$ an almost holomorphic extension of \( u_h \) of $\mathcal{G}^s$ class.

\textbf{Proof}. Our proof follows closely \cite{HelmSjoberg1975}, we adapt their arguments and deform \( \mathbb{R}^{2n} \) into $2n$ dimensional integration paths in \( \mathbb{C}^{2n} \) using Stokes' formula.

From the lemma 2.1 of \cite{HelmSjoberg1975} one has

\[
\text{Im} a(y, Z y) \geq \frac{1}{C} |\text{Im} Z y|^2 \tag{5.21}
\]

We set

\[
h(y, Z) = a(y, Z + Z y) - a(y, Z y)
\]

\( h \) is defined close to \( (y_0, 0) \) and is modulo $a$ small error a quadratic form in \( Z \). One has

\[
\partial_Z h(y,0) = 0, \quad \partial_{\overline{Z}} h(y, Z) = \mathcal{O}(1)e^{-\frac{1}{C} \text{Im}(Z + Zy)^{-1/(s-1)}} \tag{5.22}
\]
We next set

\[
R(y, Z) = 2 \int_0^1 (1 - \theta) \, h''_{Z,Z}(y,\theta Z) \, d\theta \tag{5.23}
\]

since \( t \rightarrow \exp\left(-\frac{1}{C} t^{-1/(s-1)}\right) \) is convex and increasing for small \( t \), one has

\[
\partial_{\bar{Z}} R(y, Z) = \mathcal{O}(1) \exp\left(-\frac{1}{C} \left( |\text{Im} Z |\right)^{-1/(s-1)} + \left(| \text{Im} Z_y| \right)^{-1/(s-1)} \right)
\]

We write

\[
h(y, Z) = (R(y, Z)\cdot Z,Z )/2 + \rho (y, Z) \tag{5.24}
\]

\[
|\partial_{\bar{Z}} R(y, Z)| + | \rho (y, Z)| = \mathcal{O}(1) \exp \left(- \frac{1}{C} (|\text{Im} Z|)^{-1/(s-1)} + (|\text{Im} Z_y|)^{-1/(s-1)}\right)
\]

\[
R(y, Z) = i \,^t Q(y, Z) Q(y, Z) \quad \text{with, if} \quad  Q(y_0, 0) = A, \quad ^tA R(y_0, 0) A = i I
\]

Setting

\[
\tilde{Z}(y, Z) = Q(y, Z - Z_y) \cdot (Z - Z_y)
\]

with \( Q(y_0, 0) \in GL(2n, \mathbb{C}), \, (y, Z) \rightarrow (y, \tilde{Z}(y, Z)) \) is a change of coordinates, with these coordinates one has

\[
a(y, Z) = a(y, Z y) + i/2 \langle \tilde{Z},\tilde{Z} \rangle + \rho(y, Z) \tag{5.25}
\]

One has

\[
\langle \tilde{Z}, \tilde{Z} \rangle = |\tilde{X}|^2 - |\tilde{Y}|^2 + 2i \tilde{X} \cdot  \tilde{Y}, \quad \tilde{Z} = \tilde{X} + i\tilde{Y}, \quad \rho \text{ is small}
\]

\[
| \rho(y, Z) | = \mathcal{O}(1) e^{\left( - \frac{1}{C} |\text{Im} (Z - Z y)|^{-1/(s-1)} \right)}
\]

The map \( Z \to \tilde{Z}(Z) \) has an inverse \( \tilde{Z} \to Z(\tilde{Z}) \) (in fact these maps depend on \( y \)).

We define for \( \sigma \in [0, 1] \) the integration paths

\[
\Gamma_{y, \sigma} \tilde{X} \to {Z}(\tilde{Z_\sigma}), \quad \tilde{Z}_\sigma = \tilde{X} + i \sigma g(y, \tilde{X}), \tilde{X} \in \mathbb{R}^{2n}\tag{5.26}
\]

with \( g(y, \tilde{X}) \) smooth and such that \( \mathbb{R}^{2n} \) is given by \( \tilde{Y} = g(y, \tilde{X}) \) in the new coordinates.  Using Stokes formula and deforming $\Gamma_{y,1} = \mathbb{R}^{2n} \) into $\Gamma_{y,0}$, one has

\begin{equation}
\begin{aligned}
    \text{Im} \, a(y, Z (\tilde{Z}_\sigma)) &\geq \frac{1}{C} (1 - \sigma) (|\text{Im} Z y|^2 + |\tilde{X}|^2) \\
    &\geq \frac{1}{C'} (|\text{Im} Z (\tilde{Z}\sigma)|^2) \text{  for } \sigma \in [0,1], y \in \text{neigh}(y_0, \mathbb{R}^{n}), \tilde{X} \in \mathbb{R}^{2n}
\end{aligned} \tag{5.27}
\end{equation}

Defining the \( (2n, 0) \) form in \( \mathbb{C}^{2n} \),

\begin{equation}
\begin{aligned}
    \omega_h &= f_h(Z) \, dZ_1 \, \wedge \dots \wedge\, dZ_{2n}, \quad f_h(Z) = e^{i  a(y, Z)/h} u_h(Z), \\
    &\text{since } \partial_{\bar{Z}} f_h = e^{i a(y, Z)/h} \big( \partial_{\bar{Z}} u_h + i \frac{u_h}{h} \partial_{\bar{Z}} a \big).
\end{aligned} \tag{5.28}
\end{equation}
We have

\[
\int_{\mathbb{R}^{2n}} e^{i a(y, X)/h} u_h(X) \, dX = \int_{\Gamma_{y,1}} e^{i a(y, Z)/h} u_h(Z) \, dZ_1 \wedge \dots \wedge dZ_{2n} 
\]

\[
\int_{\Gamma_{y,0}} e^{i a(y,Z)/h} u_n(Z)dZ_{1}\wedge \dots \wedge dZ_{2n} = \int_U e^{i a(y,Z_y)/h} J(\tilde{X}) u_n (Z(\tilde{X}))  e^{\frac{-1}{2h}|\tilde{X}|^2} e^{\frac{i}{h}\rho} d\tilde{X} \tag{5.29}
\]

with \( U \) a neigh of 0 in \( \mathbb{R}^{2n} \), \( J\tilde{(X)} = |\frac{\partial{Z}}{\partial{\tilde{X}}}| \),  is a Jacobian.

Using the arguments of \cite{HelmSjoberg1975} lemma 2.5, we reduce the computation of the asymptotic of the RHS of (5.29) to the case \( \rho \equiv  0 \). We pay it by a loss of a \( O_{2s-1}(h^\infty)\) remainder (see \cite{LascarLascar1997}).

Now
\[
\int_{\mathbb{R}^{2n}} e^{-\frac{1}{2h}|\tilde{X}|^2} u_n ({Z}\tilde{(X)}) J\tilde{(X)} d\tilde{X}
\]
has an asymptotic expansion. Indeed if, $(A Y, Y)$ is a non degenerate quadratic form in $\mathbb{R}^{2n}$ with $\operatorname{Im}A \geq 0$.

if $\chi(Y)$ is a $G^s_0$ truncation close 0 in $\mathbb{R}^N$. We may write,
\[
h^{-\frac{N}{2}} \int_{\mathbb{R}^N}e^{\frac{i(AY,Y)}{2h}} u_h(X + Y) dY = {\ell}_\chi(X) + r_\chi(X) \tag{5.30}
\]
\[
\ell_\chi(X) = h^{-\frac{N}{2}} \int_{\mathbb{R}^N} e^{\frac{i(AY,Y)}{2h}} u_h(X + Y) \chi(Y)dY \tag{5.31}
\]
\[
r_\chi(X) = h^{-\frac{N}{2}} \int_{\mathbb{R}^N} e^{\frac{i(A Y, Y)}{2h}} u_h(X + Y) (1 - \chi(Y)) dY.
\]

We have $r_\chi(X)$ is a $G^s$ remainder i.e. $
\lvert \partial^\alpha r_\chi \rvert \leq C^{(1 + \lvert \alpha \rvert)} \alpha!^{s}e^{(-\frac{1}{C}h^{-1/s})}$
and $\ell_\chi(X)$ is a $G^s$ symbol of zero order having the expansion
\[
\lvert \partial^\alpha ( \ell_\chi - \sum_{k=0}^{K-1} \frac{h^k}{k!} C_A \left( \frac{A^{-1}D, D}{2i} \right)^k u_h(X)) \rvert
\leq C^{ (1 + \lvert \alpha \rvert +K +1)} \alpha!^s K!^{2s-1}h^k 
\]
with $C_A = \left( \det A / 2i\pi \right)^{-1/2}$.

As a consequence,
\[
v_h(Y) = e^{\frac{-i}{h} a(y, Z_y)} \int_{\mathbb{R}^{2n}} e^{\frac{i}{h} a(y, X)} u_h(X) dX \tag{5.32}
\]

is a $G^s$ symbol of order $-n$, admitting an expansion
\[
\sum_{v=0}^\infty h^{v+n} C_{\nu,y}(D) u_h(Z_y), \tag{5.33}
\]
with $C_{\nu,y}(D)$ a differential operator of degree $\leq 2\nu$. The first term can be rewritten as
$h^n C_0(y) u_h(Z_y),$ with $ C_0(y)(2\pi)^{-n}$ a suitable branch of the square root of $\det \left(\frac{1}{i} a''_{Z, Z} (y, Z_y)\right)^{-1},$
Moreover, the last term coming from Stokes' formula is a small $\mathcal{G}^{2s-1}$ remainder, 
at least when $a(y, \cdot)$ depends analytically on $y$, which is the case here since
\[
a(y, X) = (y - z)\eta + \varphi(z), \quad X = (z, \eta). \tag{5.34}
\]

Writing this term,
\[
R(y) = \iint \limits_{U_{\sigma\in[0,1]}\Gamma_{y,\sigma}}  e^{\frac{i}{h} a(y,\eta)} \left( \partial_{\bar{Z}}{u_n} + \frac{u_n}{2} \partial_{\bar{z}}{a} \right) \wedge dZ_n \cdots \wedge dZ_{2n} \tag{5.35}
\]
And computing $\partial^{\lambda}_y R(y)$, one has
\begin{equation}
\partial^{\lambda}_y R(y) =
\sum\limits_{\substack{\lambda = \mu + \nu \\ \nu_1+ \dots + \nu_\ell = \nu \\ \ell = 1, \dots, |\nu|}} 
\frac{\lambda!}{\mu! \nu!}
\iint\limits_{U_{\sigma\in[0,1]}\Gamma_{y,\sigma}} 
e^{ i a/h} \frac{h^\ell \nu!}{\ell! \nu_1! \cdots \nu_\ell!} 
\partial^{\nu_1}_y a \cdots \partial_y^{\nu_\ell} a \partial_y^{\mu} g_h \,\wedge dZ_1 \cdots \wedge d Z_{2n}
\tag{5.36}
\end{equation}

with $ g_h = \partial_{\bar{z}} u_h + \frac{i}{h} u_h \partial_{\bar{z}} a$, for $Z \in \bigcup_{\sigma \in [0, 1]} \Gamma_{y, \sigma}$,  $\operatorname{Im} a(y, Z) \geq \frac{1}{C} \lvert \operatorname{Im} Z \rvert^2$ by (5.27), and that $\partial_{\bar{z}} u_h$ and $\partial_{\bar{z}} a$ are bounded by  $C \exp(-\frac{1}{C} |Im Z|^{-1/s-1})$ like their $\partial_y$  derivatives by design. We deduce estimates for the RHS of (5.36) since

\[
\left\lvert \frac{\partial^{\nu_1}_ya}{{\nu_1!}}\right\rvert  ...  \left\lvert 
 \frac{\partial^{\nu_\ell}_ya}{{\nu_\ell!}} \right\rvert  
\leq C^{|\nu|+1}, \quad 
h\left\lvert {\partial_y^\mu g_{h}}\right\rvert \leq C^{1+\lvert \mu \rvert} \mu!^s e^{\left(-\frac{1}{C} \lvert \operatorname{Im} Z \rvert^{\frac{-1}{s-1}}\right)}, 
\left\lvert e^{\frac{ia(y, Z)}{h}} \right\rvert \leq e^{-\frac{1}{C} \lvert \operatorname{Im} z \rvert^2}.\tag{5.37}
\]

On the other hand,
\[
\sum_{\lambda = \mu + \nu} \frac{\lambda!\nu!}{\mu! \nu! e!} h^{-\ell} M! N! h^N \lvert \operatorname{Im} Z \rvert^{-2 N} \lvert \operatorname{Im} Z \rvert^{\frac{M}{s-1}}\mu!^{s} \tag{5.38}
\]

For all \( N, M > 0 \), choosing \( N = \ell+{M'} \), with $M'>0$ given \( M = 2N(s - 1) \), one has the bound for (5.38):  
\[
C^{1 + |\lambda|} \mu!^s \nu! \ell!^{2s - 2} \leq C^{1 + |\lambda|} \mu!^{s} \nu!^{2s - 1}\leq   C^{1 + |\lambda|} \lambda!^{2s - 1}
\]
since  
\[
\frac{M! N!}{\ell!} \leq C^{1 + M'} M'!^{2s - 1} \ell!^{ 2s - 2}
\]
so \( \partial_y^\lambda R(y) \) satisfies  
\[
|\partial_y^\lambda R(y)| \leq C^{1 + |\lambda|} \lambda!^{2s - 1} e^{\left( -\frac{1}{C} h^{\frac{-1}{2s - 1}} \right)}
\tag{5.39}\]

Hence $R$ is a small $G^{2s-1}$ remainder. Finally, we have written
\[
\iint_{\mathbb{R}^{2n}} e^{\frac{i}{h} a(y, X)} u_h(X) \, dX = e^{\frac{i}{h} a(y, Z_y)} v_h + R(y), \tag{5.40}
\]
with $v_h$ on a $\mathcal{G}^s$ symbol and $R(y)$ a small $\mathcal{G}^{2s-1}$ remainder.

This ends Lemma 5.1 proof. We have the remark:
\subsection*{{{Remarks 5.2}}}
1) For obtaining the asymptotic $(5.15)$ modulo a $\mathcal{G}^{2s-1}$ remainder, we use Lemma $(5.1)$ and a reduction to the case where $q(y, \theta)$ has compact support in $\theta$. It is performed by using a cutoff and integration by parts in $z$. We write
\[
Q(a e^{\frac{i\varphi}{h}})(y) = \iint_{\mathbb{R}^{2n}} e^{\frac{1}{h}((y-z)\theta+ \varphi(z))}
 \frac{q(y, \theta)a(z)}{(2\pi h)^n} \, dz d\theta \, . \tag{5.41}
\]

$\varphi$ is a $\mathcal{G}^s$ symbol of order 0, so one has $\lvert \varphi'_y \rvert \leq C/2$ and for $\lvert \theta \rvert \geq C > 0$, 
the phase $H(y, z, \theta) = (y - z) \theta + \varphi(z)$ has no critical point in $z$ so we may perform integration by parts with 
\[
L = h \frac{H'_z}{i \lvert H'_z \rvert^2} \partial_z,
\]
splitting the integrand $q(y, 0) a(z)$ into two terms $\chi(\theta) q(y, \theta) a(z)$ and $(1-\chi)(\theta) g(y, \theta) a(z)$. 
We integrate by parts on the second term and obtain, applying the above lemma:
\[
Q(a e^{\frac{i}{h} \varphi}) = Q_\varphi^\#(a) e^{\frac{i}{h} \varphi} + R(a). \tag{5.42}
\]

$Q_\varphi^\#(a)$ is a $\mathcal{G}^s$ symbol, $R(a)$ is a $\mathcal{G}^{2s-1}$ remainder. Moreover,
\[
Q_\varphi^\#(a) \text{ has the asymptotic expansion } \sum_{\alpha \geq 0} \frac{h^{\lvert \alpha \rvert}}{\alpha!} D_\theta^\alpha \tilde{q}(y,\varphi_y') D_z (e^{\frac{i}{h}\rho}a) \rvert_{y = z},\]
(in the smooth sense). $Q_\varphi^\#(a)$ is a $\mathcal{G}^s$ symbol of order 0. 

Let us end this remark by observing that the above operators $C_{\nu, y}(D)$ of (5.33) may be computed in using $Q(y, h D_y)$ operators with polynomial symbols (see \cite{HelmSjoberg1975}).

Applying Lemma 5.1 and Remark 5.2, we have:
\[
P'a' = e^{-\frac{i}{h} \varphi}\{  h D_{Rex_1} - h^{m}  { }^tP(y, D_y) \} (a'e^{\frac{i}{h} \varphi}) \tag{5.43}
\]

rewrites on $\Gamma$:
\[
P'a' = \frac{h}{i} \partial_{Rex_1} a' + \frac{h}{i} \frac{\partial p}{\partial \eta} (y, - \varphi'_y) \partial_ya' + c'a' + R'a'. \tag{5.44}
\]

In $(5.44)$, $c'$ is a symbol of order $-1$, $R'$ a linear map $S^{m'}_s \to S_s^{m'-2}$.

We have written in $(5.44)$ the transport equation on $\Gamma$ we have to solve the eikonal equation. We set $V'$ for the real vector field:
\[
V' = \partial_{Re(x_1)} + \frac{\partial p}{\partial \eta} (y - \varphi'_z) \partial_y  + \overline{\frac{\partial p}{\partial \eta} (y - \varphi'_y)} \partial_{\bar{y}}\tag{5.45}\]

Notice we checked above that $V'|_\Gamma$ is tangent to $\Gamma$ (see $(5.11)$).

We solve inductively the equations:
\[
\frac{h}{i} V' a_0 + c a_0 = h b,\quad \frac{h}{i} V' a_j + c a_j + R a_{j-1} = 0 \quad \text{for } j \geq 1. \tag{5.46}
\]

We will construct $a'_j \in S^{-j}_s(\mathbb{C}^n\times \mathbb{C}^n)$, almost holomorphic, with respect to $\Gamma$ such that the restrictions $a_j = a'_j|_\Gamma$, satisfies $(5.46)$. We have noted $R a_{j-1}$. for $R' a'_{j-1}|_\Gamma$.

In fact, we have not used that $P(y, D_y)$ is a differential operator, only 
that $h^m P(y, D_y)$ is a zero-order Gevrey PDO on $\mathbb{R}^n$. Using the Egorov property obtained in Theorem. 2, we reduce $P(y, hD_y)$ to $P'(y, hD_y)$ having the same principal symbol such that
\[
P' = P'_0(y, hD_y) + Q', P'_0=op_h(p)
\]
with $Q'$ a $h$ PDO of order $m_0'$ small, $m_0' \leq -n - 1$. It is obtained by using the Gevrey calculus of PDO. Reproducing the same induction process as in Lemma $3.2$, we obtain  a sequence $(a'_j)_{j \geq 0} \in S^{-j}_s$, which is a formal symbol in  $\mathbb{C}^n \times \mathbb{C}^n$, close to $(x_0,y_0)$, if $a' \sim \sum_{j \geq 0} a'_j$ and  $b' = \frac{1}{h} P' a'$. One has
\[
h b' = \{ \varphi'_{x_1} - p(y, -\varphi'_y) \} a' + \{ \frac{h}{2} V'  + c ' + R'\} a'.
\]

Since on $\Gamma$ (for all $N> 0$),
\[
hb'= \frac{1}{h} \{S'_\chi (a') + R' a'_{N-1} 
+ (\frac{h}{i} V' + c' + R')
(a' - \sum_{j < N} a'_j)\}.
\]

As we have satisfied on $\Gamma$ the eikonal and transport equations. 

Taking coordinates in $\mathbb{C}^n \times \mathbb{C}^n$ such that $\Gamma$ is $\{y' = 0\}$, with $z' = x' + i y'$, with $x', y' \in \mathbb{R}^{2n}$, one has for all $\alpha $ and , $N > 0$:
\[
\lvert \partial_{x'}^\alpha b' \rvert \leq C^{1 + \lvert \alpha \rvert + N} N!^s   \alpha!^s h^N.
\]

Indeed, $S'_\chi(a')$ is a small $G^s$ remainder, $\rho'^N = a' - \sum_{j < N} a'_j$, is a $-N$ symbol, $R'a'_{N-1}$ is a $-N-1$ symbol and $so$ 
\[
\left(\frac{h}{i} V' + c' + R'\right)(\rho'^{N}).
\]

Moreover, we note too that the Boutet de Monvel \cite{BoutetDeMonvel1967} proof of the existence of solutions to Borel problems in Gevrey calculus allows us to choose $a'(x, y, h)$ almost holomorphic with respect to $\Gamma$. Using this fact and the above remarks, we write:

\[ h b'(x, y, h) = e^{-\frac{i}{h} \varphi} \{ ( hD_{Re(x_1)}) - h^{m} {}^t P(y, D_y))\} (a'e^{\frac{i}{h} \varphi}) \tag{5.47} \]

$b'$ satisfies in coordinates $z' = x' + i y'$:
\[
b'(x' + i y') = \sum_{\lvert \alpha \rvert \leq n} \frac{\partial^\alpha_{x'} b'(x')(iy')^\alpha}{\alpha!}  + R_N(x', y') \tag{5.48}
\]

with
\[
R_N(x', y') = \sum_{|\alpha| \leq N} \frac{N!}{\alpha!} \int_{0}^{1}  \frac{(1-s)^{N-1}}{(N-1)!} \partial y'^{\alpha}b' (x' + iy') ds {y'^\alpha}
\]

Each $b'^{\alpha}(x')$ satisfies for all $M > 0$,

\[
\lvert b^{(\alpha)}(x') \rvert \leq C^{1 + \alpha + M} h^{M/s} \alpha!^{s}  M!
\]

Since 
\[
\lvert b^{(\alpha)}(x') \rvert \leq C^{1 + \lvert \alpha \rvert} \alpha!^s e^{\left(-\frac{1}{C} h^{-1/s}\right)},
\]
from $(5.48)$ one has for $\epsilon > 0$ given, for $\lvert y' \rvert \leq C h^\epsilon$, estimates for $\lvert \alpha \rvert \leq N$:
\[
\frac{\lvert b^{(\alpha)}(x')(y')^\alpha \rvert}{\alpha!} 
\leq C^{M + N + 1} h^{ M/s} \alpha!^{s-1} h^{\epsilon\lvert \alpha \rvert }M! \tag{5.49}
\]

Choosing $\epsilon = \frac{s-1}{s}$, $M = \epsilon s (N-\lvert \alpha \rvert)$, and defining:
\[
F_N = \sum_{\lvert \alpha \rvert \leq N} \frac{\partial_x^\alpha b'(x')(\epsilon y')^\alpha}{\alpha!},
\]
we have:
\[
\lvert F_N \rvert \leq C^{1 + N} N!^{s-1} h^{\epsilon N}.
\]
Since $R_N$ satisfies:
\[
\lvert R_N \rvert \leq C^{N+1}N!^{s-1} \lvert y' \rvert^N \leq C^{N+1} N!^{s-1} h^{\epsilon N},
\]
so for $\lvert y' \rvert \leq C h^\epsilon$, one has:
\[
\lvert b'(x' + i y') \rvert \leq C^{N+1} N!^{s-1} h^{\epsilon N}. \tag{5.50}
\]

Taking in (5.50) the infimum over $N$, one has:
\[
\lvert b'(x' + i y') \rvert \leq C' \exp\left(-\frac{1}{C'}  h ^{-1/s}\right) \text{ for } \lvert y' \rvert \leq C h^\epsilon. \tag{5.51}
\]

Now for $\lvert y' \rvert \geq C h^\epsilon$, using similar arguments, one has:
\[
\lvert b'(x' + i y') \rvert \leq C' \exp\left(-\frac{1}{C} \lvert y' \rvert^{\frac{-1}{s-1}}\right).
\]

One has in all cases
\[
\lvert b'(x' + i y') \rvert \leq C \{ e^{(\frac{-1}{C} h^{-\frac{1}{s}})} + e^{(-\frac{1}{C} \lvert y' \rvert^{\frac{-1}{s-1}})} \} \tag{5.53}
\]

One deduces, if $u \in \mathcal{D'}(\mathbb{R}^n)$ is independent of $h \in (0, 1]$, estimates:
\[\lvert \partial_{x}^\gamma T_{b'}u(x) \rvert \leq  C_\gamma e^{\frac{1}{h}\phi(x)} e^{(-\frac{1}{C} {h}^{\frac{-1}{2s-1}})}, \quad x \in \text{neigh}(x_0, \mathbb{C}^n), \tag{5.54}\]
with $\phi(x) = \sup_{y \in \text{neigh}(y_0, \mathbb{R}^n)} -\operatorname{Im} \varphi(x, y)$. (In fact, $\phi(x) = -\operatorname{Im} \varphi(x, y(x))$.)

We have Theorem 3's statement since $(5.54)$ is the desired result if one has adjusted $\phi = \phi_0$ in choosing the convenient primitive of $\omega = \xi dx - \eta dy |_\Gamma$.

A corollary could be the Gevrey $(2s-1)$ propagation of singularities. We refer the reader to \cite{Hitrik2023personal}.

\appendix
\section{Gevrey WF set of Lagrangian distributions}

Let $u(x, h)$ be 
\[
u(x, h) = \int e^{\frac{i}{h} \phi(x, \theta)} a(x, \theta, h) \, d\theta, \tag{A.1.1}
\]
a Lagrangian distribution on $\mathbb{R}^n$ with:

$\phi(x, \theta) \in G^s(\mathbb{R}^n \times (\mathbb{R}^N \setminus \{0\}))$, a homogeneous phase of degree $1$ with respect to $\theta$, variables with $\operatorname{Im} \phi \geq 0$, $d \phi \neq 0$, $ \cdot a(x, \theta) \in \tilde{S}^{m_0,0}_{s}(\mathbb{R}^n \times \mathbb{R}^n)$, a Gevrey symbol of order $(m_0, 0)$, i.e., 
\[
\lvert \partial_x^\alpha \partial_\theta^\beta a(x, \theta) \rvert 
\leq C^{1 + \lvert \alpha \rvert + \lvert \beta \rvert} h^{-m_0} \alpha!^s \beta!^s (1 + \lvert \theta \rvert)^{-|\beta|}, \tag{A.1.2}
\]

If $\operatorname{WF}_h(u)$ is the semiclassical wavefront set of $u$, one has at least in $\xi\neq 0$:
\[
\operatorname{WF}_h(u) \subset \{(x, \phi'_x) \mid \phi_\theta(x, \theta) = 0, (x,\theta)\in \operatorname{cone supp a} \}
\]

In the Gevrey case $(s > 1)$, one has a similar result:
\subsection*{{Lemma A.1.2}} Let $u(x, h)$ as above, one has:
\[
\operatorname{WF}_{s, h}(u) \subset \{(x, \phi'_x(x, \theta));  \phi'_\theta(x, \theta) = 0, (x, \theta) \in {cone supp a} \} = \Lambda.
\]

\textbf{Proof.} The proof follows the one of Proposition 8.1.9 of \cite{Hormander1985}, adapted to the Gevrey case.

Assume $(x_0, \xi_0) \notin \Lambda$, $\psi \in \mathcal{G}^s_0(\mathbb{R}^n)$ supported close to $x_0$. Let $V$ be a conic closed subset such that:
\[
V \cap \{\phi'_x(x, \theta); \xi \in \operatorname{supp} \psi, (x,\theta) \in \operatorname{cone suppa}, \phi'_\theta(x, \theta) = 0\} = \varnothing.
\]

We aim at proving that for any  $\xi $ $\in V(\xi_0) \subset\subset V$, one has:
\[
\widehat{\psi u}(\xi/h) = \mathcal{O}(1)e^{-\frac{1}{C} h^{-1/s}}, \quad \text{for } \xi \in V(\xi_0). \tag{A.1.3}
\]

Moreover:
\[
\widehat{\psi u}(\xi/h) = \int e^{\frac{i}{h} (\phi(x, \theta) - x \xi)} a(x, \theta) \psi(x) \, dx \, d\theta. \tag{A.1.4}
\]
We perform in the integrand (A.1.4) a $\mathcal{G}^s$ dyadic partition of unity in $\theta \in \mathbb{R}^N$. It has been already used in previous sections.

We rewrite (A.1.4):
\[
\widehat{\psi u}(\xi/h) = \sum_{\nu=0}^\infty \int e^{\frac{i}{h} (\phi(x, \theta) - x \xi)} a(x, \theta) \chi_{\nu} (\theta) \psi(x) \, dx \, d\theta. \tag{A.1.5}
\]

We set $R_\nu = R = {2^{\nu - 1}}$ and write for $\nu \geq 1$ the terms in $A.1.5$ as: \\
\( R^N  \int e^{\frac{i}{h} (R\phi(x, \theta) - x \xi)} \psi(x)a(x, R\theta) \chi(\theta) dx d\theta \)

We evaluate each term by Gevrey non-stationary phase lemma 2.1 of these notes.  Let $\Phi(x, \xi, \theta)$ be the phase:
\[
\Phi(x, \xi, \theta) = \frac{(R \phi(x, \theta) - x \xi)}{R + |\xi|}.  \tag{A.1.6}
\]

We have, for $\xi \in V$, $(x, \theta) \in \text{supp} \, \psi(x) a(x, R \theta) \chi(\theta)$, the bound:
\[
|\Phi'_x| + |\Phi'_\theta| \geq c'  \frac{R |\theta| + |\xi|}{R + |\xi|} \geq c > 0. \tag{A.1.7}
\]

Deducing from (A.1.7), the bound:
\[
C A^k k!^s h^{-m_0}h^kR^N (R + |\xi|)^{-k} h^{N},
\] for $k$ large enough with respect to $N$. So, one has the bound:
\[
C A^k k!^sR^{-1} h^{-m_0+k} (1+|\xi|)^{-k+1} \text{  for } k\geq k_0,\xi \in V.
\]

We may \underline{sum} these bounds since:
\[
\sum_{\nu \geq 1}^\infty 2^{-\nu+1} = 2 \quad (R = 2^{\nu-1}).
\]

Hence, we deduce for $\widehat{\psi u}(\xi/h)$ the bound $|\widehat{\psi u}(\xi/h)| \leq C (e^{-\frac{1}{C} h^{-1/s} })$ for $\xi \in V(\xi_0) \subset \subset {V}$, which is $(A.2.3)$.

We end with a remark:

\subsection*{{Remark A.1.3}} 

We may drop the homogeneity assumption of $\phi$ with respect to $\theta \in \mathbb{R}^N$ if we assume $a(x, \theta) \in S^{m_0}_s$ and has its $\theta$-support in a fixed compact of $\mathbb{R}^N$. In this case, we drop the dyadic partition of unity in the proof.

\section{Borel Lemma in Gevrey Symbol Classes}

We introduced earlier classes $\tilde{S}_s^{m'}$ ( resp $S_s^{m'})_s$. Formal symbols are sequences $(a_j)$ of symbols in $\tilde{S}_s^{m'-j}$ (resp. $S_s^{m'-j})$ such that

\[
\lvert \partial_x^\alpha \partial_\xi^\beta a_j \rvert \leq C^{1+\alpha + \beta + j} \alpha!^s \beta!^s j!^sh^{-m'+j}\langle \xi \rangle^{-j-|\beta|} \tag{(A.2.1) i) }
\]
\[
\leq C^{1+\alpha + \beta + j} \alpha!^s \beta!^s j!^s h^{-m'+j}. \tag{(A.2.1) ii)}
\]

Using the result of Boutet de Monvel-Krée \cite{BoutetDeMonvel1967} from the Carleson Theorem \cite{Carleson1964}, we may realize Gevrey formal symbols into symbols. It is a consequence of a Borel argument.
\subsection*{{Proposition A.2.1}}

Let $(a_j)$ be a formal symbol in $S_s^{m'}$ (resp. in $\tilde{S}_s^{m'})$.  There exists $a \in S_s^{m'} \text{resp.} \tilde{S}_s^{m'})$ such that for $N > 0$, $\alpha$, $\beta$ one has

\[
\lvert \partial_x^\alpha \partial_\xi^\beta ( a - \sum_{j < N} a_j )\rvert 
\leq C^{1 + \alpha + \beta + N} \alpha!^s \beta!^s N!^s h^{-m' + N}. \tag{A.2.2}
\]
resp.
\[
\leq C^{1 + \alpha + \beta + N} \alpha!^s \beta!^s N!^s h^{-m' + N} \langle \xi \rangle^{-|\beta| - N}. \tag{A.2.3}
\]

First we prove (A.2.2).

\underline{Proof}

The sequence $(a_j(x, \xi)){h^{-j} j!}$ is a $(s+1)$ sequence in the sense of Boutet de Monvel–Krée \cite{BoutetDeMonvel1967}. By the Carleson Theorem, we can solve the following Borel problem. We construct 

\[
g(t, x, \xi, h) \in G^{s-1}(\mathbb{R}, S_s^{m'}) \quad \text{such that for all } j \geq 0,
\]

\[
\partial_t^j g(0, x, \xi/ h) = j! h^{-j} a_j(x, \xi, h). \tag{A.2.4}
\]

Setting $a(x, \xi, h) = g(h, x, \xi, h)$ for small $h \in (0, 1]$, one has by Taylor's formula:

\[
\lvert \partial_x^\alpha \partial_\xi^\beta (a - \sum_{j < N} a_j) \rvert 
\leq \sup_{\theta \in (0, 1)} \lvert \partial_t^N \partial_x^\alpha \partial_\xi^\beta g(\theta h, x, \xi, h) \rvert 
\frac{h^N}{N!}. \tag{A.2.5}
\]

From (A.2.5), one has

\[
\lvert \partial_x^\alpha \partial_\xi^\beta (a - \sum_{j < N} a_j) \rvert 
\leq C^{1+\alpha+\beta+N} N!^s \alpha!^s \beta!^s h^{-m' + N}. \tag{A.2.6}
\]

g being Gevrey $s+1$ in $t, s$ w.r.t. $(x,\xi)$ one has by (A.2.5)(A.2.6).  For obtaining a similar result for the class $\widetilde{S}_s^{m'}$, we use this result and a dyadic partition in $\xi$ (see \cite{Hormander1985}). Let us write with $\chi_0 \in \mathcal{G}_0^s(\mathbb{R}^n)$, $\chi \in \mathcal{G}_0^s(\mathbb{R}^n \setminus 0)$:

\[
1 = \chi_0(\xi) + \sum_{\nu > 0} \chi(2^{-\nu} \xi), \tag{A.2.7}
\]

and set 

\[
a_0^j(x, \xi) = \chi_0(\xi) a^j(x, \xi), \quad 
a_v^j(x, \xi) = \chi(2^{-\nu} \xi) a^j(x, \xi).
\]

The sequence $(a^j)_{j \geq 0}$ being a formal symbol in $\widetilde{S}_s^{m'}(\mathbb{R}^m \times \mathbb{R}^m)$, we set then:

\[
\tilde{a}_\nu^j(x, \xi) = a_{\nu}^j(x, 2^\nu \xi), \text{ one has }\]

\[ \lvert \partial_x^\alpha \partial_\xi^\beta \tilde{a}_v^j \rvert \leq C^{1 +\alpha + \beta + j} 
\left( \frac{h}{2^\nu} \right)^j \alpha!^s \beta!^s j!^s h^{-m'} \tag{A.2.8}
\]

By (A.2.8), one has that the sequence $(\tilde{a}_v^j)_{j \geq 0}$ is a formal ${S}_s$ symbol in ${S}_s^{m'}$.

(See (A.2.1)) with the small parameter $\widetilde{h} = \frac{h}{2^\nu}$. The $a_\nu^j$ have their support  on a fixed compact where $|{\xi}| \sim 1$. 

We realize the sequence $(a_\nu^j)_{j \geq 0}$ into a symbol $b_\nu \in \widetilde{S}_s^0$ with $\xi$  support compact in $[{\xi}] \sim 1$.  We set now $a_\nu(x, \xi) = b_\nu(x, \frac{\xi}{2^\nu})$ 
and one has

\[
\left| \partial_x^\alpha \partial_{{\xi}}^\beta 
\left( a_\nu - \sum_{j <N} a_\nu^j \right) \right| 
\leq C^{\alpha + \beta + N+1} 
\left( \frac{h}{2^\nu} \right)^N \langle \xi \rangle^{-|\beta|} \alpha!^s \beta!^s N!^s h^{-m'}, 
\]

uniformly in $\nu, N, h$. Setting $R_\nu^N = a_\nu - \sum_{j < N} a_\nu^j$, the sequence $(R_\nu^N)_{\nu \geq 0}$ may be summed in $\Sigma^{-N,-N+\frac{1}{2}}_s$ (we have put $m' = 0$) and set $R^N = \sum_\nu R_\nu^N$.

We realize too $(a_0^j(x, \xi))_{j \geq 0}$ into $a_0(x, \xi)$ with $\xi$ compact support. Finally we set $a = \sum_\nu a_\nu + a_0$.

Let us check that for some $C > 0$, one has for all $N > 0$, $\alpha, \beta$,
\[
\left| \partial_x^\alpha \partial_\xi^\beta \left( a - \sum_{j < N} a^j \right) \right| 
\leq C^{\alpha + \beta + N + 1} \alpha!^s \beta!^s N!^s 
\langle \xi \rangle^{-|\beta| - N + \frac{1}{2}} h^N \tag{A.2.9}
\]

Writing $a - \sum_{j < N} a^j = \sum_\nu \left( a_\nu - \sum_{j < N} a_\nu^j \right) + R_0^N 
= \sum_\nu R_\nu^N + R_0^N = R^N + R_0^N$,

one has
\[
\left| \partial_x^\alpha \partial_\xi^\beta 
\left( a - \sum_{j < N} a^j \right) \right| 
\leq C^{\alpha + \beta + N + 1} \alpha!^s \beta!^s N!^s h^N 
\langle \xi \rangle^{-|\beta| - N + \frac{1}{2}} \tag{A.2.10}.
\]

Now, moving \( N\) into \( N + 1 \), one has since \( R^N = a^N + R^{N+1} \), estimates for \( R^N \),
\[
\left| \partial_x^\alpha \partial_\xi^\beta R^N \right| \leq C''^{|\alpha|+|\beta|+N+1} (N+1)!^s h^N \alpha!^s \beta!^s \langle \xi \rangle^{-|\beta|-N} \tag{A.2.11}.
\]

Writing \((N+1)!^s \leq C'''^{N+1} N!^s\) for  \( C^{'''} > 0 \), choosing \( C \geq \sup(C', C'', C''') \), 

one has
\[
\left| \partial_x^\alpha \partial_\xi^\beta \left( a - \sum_{\delta < N} a^j \right) \right| \leq C^{\alpha+\beta+N+1} \alpha!^s \beta!^s N!^s h^N \langle \xi \rangle^{-|\beta|-N} \tag{A.2.10}.
\]

So we have proved the Borel Lemma (A.2.1).\\
\\

\underline{Disclosure statement}

The authors do not work for, consult, own shares in or receive funding from any company or organization that would benefit from this article, and have disclosed no relevant affiliations beyond their academic appointment.

\bibliographystyle{unsrt} 

\bibliography{bibliography}

\end{document}